\documentclass[oneside, 12pt]{amsart}
\usepackage{amscd, amssymb, amsmath, mathrsfs}
\usepackage[english]{babel}
\usepackage{booktabs}
\usepackage{tikz-cd}
\usepackage{url}
\usepackage[pdftex, colorlinks=true,  citecolor=blue, linkcolor=blue, linktocpage=true]{hyperref} 

\usepackage{booktabs,longtable}
\usepackage{multirow}

\usepackage{multirow}
\usepackage{enumerate}
\usepackage{tikz}
\usetikzlibrary{cd,arrows,positioning} 
\tikzset{>=stealth}
\tikzcdset{arrow style=tikz}
\tikzset{link/.style={column sep=1.8cm,row sep=0.16cm}}
\tikzset{map/.style={row sep=0em, column sep=0em}}
\usepackage{amssymb}
\usepackage{array}
\usepackage{mathrsfs}
\usepackage{xspace}
\usepackage{MnSymbol}
\usepackage{mathtools}
\usepackage{color}
 \usepackage[normalem]{ulem} 
\usepackage[curve,all]{xy}
\xyoption{matrix}
 \xyoption{curve}
 \xyoption{color}
 \xyoption{line}
\xyoption{arc}

\usepackage{soul}

\usepackage[shortlabels]{enumitem}
\setlist{wide}
\setlist[enumerate]{label=\rm{(\arabic*)}}
\setlist[enumerate,2]{label=\rm({\it\roman*})}
\setlist[itemize]{label=\raisebox{0.25ex}{\tiny$\bullet$}}

\newtheorem{theorem}{Theorem}[section]

\newtheorem{lemma}[theorem]{Lemma}
\newtheorem{proposition}[theorem]{Proposition}
\newtheorem{corollary}[theorem]{Corollary}
\newtheorem*{mainthm*}{Theorem}
\newtheorem{question}[theorem]{Question}

\theoremstyle{definition}
\newtheorem{definition}[theorem]{Definition}

\newtheorem{examples}[theorem]{Examples}

\newtheorem{remark}[theorem]{Remark}

\newtheorem*{acknowledgement}{Acknowledgement}

\theoremstyle{remark}

\newcommand{\incl}[1][r]{\ar@<-0.2pc>@{^(-}[#1] \ar@<+0.2pc>@{-}[#1]}

\renewcommand{\P}{\mathbb{P}}

\newcommand{\Spec}{\mathrm{Spec}}

\newcommand{\p}{\mathbb{P}}

\newcommand{\Q}{\mathbb{Q}}

\newcommand{\A}{\mathbb{A}}

\newcommand{\C}{\mathbb{C}}
\newcommand{\F}{\mathbb{F}}

\newcommand{\Z}{\mathbb{Z}}
\newcommand{\R}{\mathbb{R}}

\renewcommand{\O}{\mathcal{O}}

\newcommand{\Gr}{\mathrm{Gr}}
\newcommand{\LGr}{\mathrm{LGr}}

\newcommand{\xr}{X(\mathbb{R})}

\DeclareMathOperator{\Fa}{F}

\DeclareMathOperator{\rk}{rk}

\DeclareMathOperator{\defo}{def}
\DeclareMathOperator{\Bl}{Bl}

\DeclareMathOperator{\Pic}{Pic}

\DeclareMathOperator{\Gal}{Gal}

\newcommand{\abs}[1]{\left\lvert#1\right\rvert}

\title[A Rationality Criterion for Real Fano Threefolds]{
A Rationality Criterion for Real Fano Threefolds}
\date{\today}

\author[Andrea Fanelli]{Andrea Fanelli}
\address{Univ. Bordeaux, CNRS, Bordeaux INP, IMB, UMR 5251, F-33400 Talence, France.}
\curraddr{}
\email{andrea.fanelli@math.u-bordeaux.fr}

\author[Fr\'ed\'eric Mangolte]{Fr\'ed\'eric Mangolte}
\address{Aix Marseille University, CNRS, I2M, Marseille, France}
\email{frederic.mangolte@univ-amu.fr}

\let\origmaketitle\maketitle
\def\maketitle{
  \begingroup
  \def\uppercasenonmath##1{} 
  \let\MakeUppercase\relax 
  \origmaketitle
  \endgroup
}


\begin{document}


\begin{abstract}
We study the connectedness of the real locus of smooth geometrically rational Fano threefolds and prove a sufficient criterion of $\R$-rationality.
\end{abstract}

\maketitle
\tableofcontents


\section{Introduction}
\label{Introduction}

Let $X$ denote a real projective scheme. A general and very hard problem consists in understanding connections
between the algebraic properties of $X$ and the topological properties of its real locus $X(\R)$, endowed with the euclidean topology. 

A first instance of what Kollár calls \emph{the recognition problem} in \cite{K01} is the following easy fact: let $X$ be smooth and $\R$-rational (i.e.\ birational to $\P^n$ over $\R$), then the real locus $X(\R)$ is connected.
\\The classical work of Comessatti \cite{Co12} on the real algebraic surfaces combined with the MMP for real surfaces shows that, in dimension 2, connectedness of the real locus characterizes rational varieties among geometrically rational ones (i.e.\ such that $X_\C$ is birational to $\P^n$ over $\C$). In \cite[III.4]{S89}, the following result appears.

\begin{mainthm*}[Comessatti]
Let $X$ be a smooth real projective surface. Suppose that $X$ is geometrically rational. Then the following are equivalent:
\begin{itemize}
\item[$\bullet$] $X$ is $\R$-rational;
\item[$\bullet$] $X(\R)$ is non-empty and connected.
\end{itemize}

\end{mainthm*}

Very little was known about the topology of the real locus for smooth real algebraic threefolds before the late 90s. The fundamental work of Koll\'ar \cite{KoI,KoII,KoIII,KoIV} on real uniruled threefolds and the minimal model program (see also \cite{K98}) and subsequent works on Koll\'ar conjectures \cite{V99,EGH,hm1,hm2,cm1,cm2,mw1} (see also \cite{M14} for an introduction to the topic) were major contributions to the field.

The literature on the explicit topology of the real locus for threefold Mori fiber spaces (i.e.\ end products of the minimal model program for uniruled varieties) focuses on conic bundles and del Pezzo fibrations and leaves aside minimal Fano threefolds. To our knowledge, the only exceptions are \cite{KS04}, \cite{Kra09} and \cite{K18}. In any case, even determining the number of connected components for the real loci of real Fano threefolds is a hard problem in general.

Returning to the rationality problem, it is known that already in dimension three, the connectedness of the real locus is not enough to guarantee rationality, and, in recent years, several works appeared to study the obstructions to $\R$-rationality of geometrically rational threefolds (see \cite{BW20,HT21, HT21b,BW23,FJ24, FJSVV2, BP24-pre, CP24}).

In this paper, we will focus on \emph{smooth Fano threefolds} and establish a sufficient criterion for $\R$-rationality.

\begin{mainthm*}[= Theorem~\ref{thm_newmain}]
Let $X$ be a smooth geometrically rational real Fano threefold with $\xr\ne\varnothing$.
If no complex deformation of $X_\C$ admits real forms whose real locus has at least two connected components, then $X$ is rational.
\end{mainthm*}
For more details on smooth geometrically rational real Fano threefold for which some deformation of $X_\C$ admit real forms whose real locus has at least two connected components, see Table~\ref{tab.recap}. 
The proof is via a case-by-case analysis,  exploiting the classification of complex families of smooth Fano threefold completed by Iskovskih  in \cite{I79} when $\rho=1$ and by 
Mori and Mukai for higher Picard rank in \cite{MM86,MM03} (see also \cite{IP99}, \cite{Mat95}, \cite{Mat-erratum} and \cite{fanography}).

For the actual statement of Theorem~\ref{thm_newmain}, we introduce an invariant $s_X$ (see Definition~\ref{def:sx}), which encodes the maximal number of connected components for the real locus of a variety $X$, up to complex deformation. We stress that, in general, it is hard to determine the maximal number of connected components of the real locus in a given deformation class. For example, this problem is still open for smooth real quintic surfaces in $\mathbb{P}^3$, see \cite{IK96,Or01}.

In the process of proving our main theorem, we reserve special care to {\it minimal} (i.e.\ with real Picard rank 1) smooth geometrically rational Fano threefolds, which play a central role in the birational classification of algebraic varieties and the minimal model program. We first survey some known results in the field and moreover produce few new examples with connected and disconnected real loci, via birational geometry and smoothing techniques (see Proposition~\ref{prop:X16}) or via explicit computations (see Propositions~\ref{prop:X33-connected} and \ref{prop:X33}).

This paper is organized as follows. In Section~\ref{sec_st} we deduce some bounds on the number of connected components for the real loci of Fano threefolds, as a consequence of the Smith-Thom theory. In Section~\ref{sec_complex} we study the real locus for minimal geometrically rational Fano threefolds, producing new examples. The main theorem is proved in Section~\ref{sec:criterion} and the final table in Secion~\ref{sec_table} summarizes the the results concerning families of Fano threefolds with many connected components.

Unless otherwise stated, all real varieties are geometrically integral and projective.
For the classification of complex families of smooth Fano threefolds, we follow the Mori--Mukai numbering of the 105 families, and write ``family \textnumero $m.n$'', where $m$ is the rank of the Picard group of the threefold, ranging from 1 to 10, and $n$ is the entry number in the list.
If $X$ is a smooth real Fano threefold, we will say that $X$ belongs to family \textnumero $m.n$ if $X_{\mathbb{C}}$ does. You can go to the webpage~\cite{fanography} to easily go through the classification.

\begin{acknowledgement}
We thank Hamid Abban, Olivier Benoist, Erwan Brugall\'e, Fabrizio Catanese, Ivan Cheltsov, Jules Chenal, Ilia Itenberg, Lena Ji, Viatcheslav Kharlamov, János Kollár, Andrea Petracci, Antoine Pinardin for fruitful discussions.
We thank Zhijia Zhang for creating an explicit example in Proposition~\ref{prop:X33}    and providing the picture.
The first author is currently supported by the ANR project ``FRACASSO'' ANR-22-CE40-0009-01.
\end{acknowledgement}


\section{
Upper bounds on the number of connected components}
\label{sec_st}

In this section, we deduce from the classical Smith-Thom theory some upper bounds for the number of connected components for the real loci.

Let $X$ be a real smooth variety of dimension $n$. First assume that $n=1$, that is $X$ is a real smooth projective curve of genus $g$. Then by a famous theorem of Harnack and Klein\footnote{
For notations and results used in this section, we refer to \cite[§3.2 and §3.6]{M17,M20}.}, the number $s$ of connected components of the real locus $X(\R)$ is at most $g+1$. Furthermore, for any $g\geq 0$, Harnack constructed a real smooth projective curve of genus $g$ with $s=g+1$. 
The Harnack-Klein inequality 
    $s\leq g+1$

was later generalized for any $n$. It's the 
\emph{Smith-Thom inequality} on Betti numbers with coefficients in $\Z/2$:
\begin{equation}
\label{eq.smiththom}
\sum_{l=0}^n b_l(\xr,\Z/2)\leqslant \sum_{k=0}^{2n} b_k(X(\C),\Z/2). 
\end{equation}

Using Galois cohomology groups, we can refine~\eqref{eq.smiththom}. 
Let $\sigma=\operatorname{id_X}\times\Spec(z\to \bar z)$ be the anti-regular involution on $X_\C=X\times \Spec\ \C$ and denote by $G=\Gal(\C/\R)=
\langle \sigma\rangle\simeq \Z/2$ the Galois group. 
Recall 
that for a $G$-module $M$, the Galois cohomology pointed sets $H^i\left(G,M\right)$, $i>0$, are in fact groups and even $\Z/2$-vector spaces. Moreover, 
$H^2\left(G,M\right)=\ker(1-\sigma)/\operatorname{Im}(1+\sigma)$ and 
we have the \emph{Borel–Swan inequality}:
\begin{equation}
\label{eq.gmz1}
\sum_{l~\text{even}}b_l(\xr,\Z/2 )
\leqslant \sum_{k=0}^{2n}\dim_{\Z/2} H^2\left(G,H_k(X(\C) ,\Z)\right) \;.
\end{equation}

When $X$ is a Fano variety two inequalities can be deduced from the previous ones. We denote by $h^{i,j}(X)=\dim_\C H^{i,j}(X_\C)$ the Hodge numbers of $X_\C$.
\begin{proposition}
\label{prop.st.gm}
Let $X$ be a real smooth Fano threefold, and $s$ be the number of connected components of the real locus $\xr$. 
Then
\begin{equation}
\label{eq.bound1}
s\leqslant 1+h^{1,2}(X)+\rho(X_\C)\;.
\end{equation}

Furthermore 
letting $\lambda_X=\rk((1+\sigma)\operatorname{Pic}(X_\C))$, we have
\begin{equation}
\label{eq.bound2}
s\leqslant 1+h^{1,2}(X)+\rho(X_\C)-2\lambda_X\;.
\end{equation}
\end{proposition}

\proof
Inequality \eqref{eq.bound1} follows directly from Smith-Thom inequality. 
Indeed the following holds (see e.g. \cite{IP99}):
\begin{enumerate}
\item The $6$-dimensional real manifold $X(\C)$ is simply-connected, and its cohomology is without 
torsion;
\item $H^i(X_\C,\O_{X_\C})=0$ for $i>0$;
\item $H^2(X(\C),\Z)\simeq \operatorname{Pic}(X_\C)$. 
\end{enumerate} 
Hence Betti numbers mod $2$ are equal to usual Betti numbers and we have
$$b_0(X(\C))=b_6(X(\C))=1,\ \  b_1(X(\C))=b_5(X(\C))=0,$$ $$b_3(X(\C))=2h^{1,2}(X),\ \  b_2(X(\C))=b_4(X(\C))=\rho(X_\C).$$ 

If non empty, the real locus $\xr$
is a $\Z/2$-oriented compact manifold of real dimension $3$ and, by Poincaré duality, one has $s=b_0(\xr)=b_3(\xr)$ and $b_1(\xr)=b_2(\xr)$. From  
\eqref{eq.smiththom}, 
we get 
$$
2s+2b_1(\xr)\leqslant 2+2\rho(X_\C)+2h^{1,2}(X)
$$
hence
$$
s+b_1(\xr)\leqslant 1+\rho(X_\C)+h^{1,2}(X)\;.
$$

We now prove \eqref{eq.bound2} using the inequality~\eqref{eq.gmz1}. 
The group $G$ acts as orientation-reversing involution on the $6$-dimensional real manifold $X(\C)$. 

Thus $H^2\left(G,H_0(X(\C) ,\Z)\right)\simeq 
\Z/2\Z$ and $H^2\left(G,H_6(X(\C) ,\Z)\right)
=0$. 
Recall that $H^3(X(\C) ;\C)=H^{2,1}(X_\C)\oplus H^{1,2}(X_\C)$ and $\sigma^*H^{2,1}(X_\C)=H^{1,2}(X_\C)$, see \cite[Proposition~D.3.17]{M17,M20}, 
then 
we have
$$
\dim_{\Z/2}H^2\left(G,H^3(X ,\Z)\right)\leqslant  h^{1,2}(X)\;.
$$

Now by definition, 
we have $\dim_{\Z/2}H^2\left(G,H_4(X(\C) ,\Z)\right)=\rho(X)-\lambda_X$. Recalling that $\Z(1)$ denotes the $G$-constant sheaf $\Z$ on which $\sigma$ acts as $m\mapsto -m$, we have by  Poincar\'e duality in group cohomology, 
$H^2\left(G,H_2(X(\C) ,\Z)\right)\simeq H^2\left(G,H_4(X(\C) ,\Z(1))\right)$. 
We deduce that $H^2\left(G,H_4(X(\C) ,\Z(1))\right)=\rho(X_\C)-\rho(X)-\lambda_X$.

Inequality \eqref{eq.gmz1} 
then gives
$$
s+b_1(\xr)\leqslant 1+h^{1,2}(X)+\rho(X_\C)-2\lambda_X\;.
$$
\endproof

\begin{remark}
In fact,  \eqref{eq.bound1} and \eqref{eq.bound2} are in general not sharp. In many cases, $s=1$ but $h^{1,2}(X)>0$. See also the table in Section~\ref{sec_table}.
\end{remark}


\section{Real loci of minimal smooth geometrically rational real Fano threefolds}
\label{sec_complex}

In this section we will study smooth real Fano threefold verifying the following two hypothesis:
\begin{equation}\label{assumptions}
\text{$X$ is geometrically rational and $\rho(X)=1$.}
\end{equation}

We refer to them as {\it minimal} smooth geometrically rational real Fano threefold and focus on them separately since play a special role, as end products of the minimal model program. 
\\We recall here some classification results, starting with a stronger assumption on our Fano variety $X$: namely, we require that $\rho(X_\C)=1$ (i.e.\ $X$ to be {\it geometrically minimal}). 

\begin{lemma}\label{lemma_real_forms}
Let $X$ be a smooth geometrically rational real Fano threefold verifying $\rho(X_\C)=1$. Then $X_\mathbb{C}$ belongs to one of the families in Table~\ref{tab1}. Moreover, the general member of the starred families are not rational.
\end{lemma}

\newcounter{tables}\setcounter{tables}{0}
\newcommand{\ltables}[1]{
\thetables\label{#1}}

\begin{center}
\phantomsection\refstepcounter{tables}
\renewcommand{\arraystretch}{1.2}
\begin{longtable}{@{}l|c|c|c|c|p{9.5cm}@{}}
\toprule
\textnumero $m.n$			& $\iota$	& $d$ & $ g$ & $h^{1,2}$	 &  Description of $X_\C$
\\
\midrule
\endfirsthead
\toprule
\textnumero $m.n$		& $\iota$	& $d$ & $ g$ & $h^{1,2}$	 &  Description of $X_\C$
 \\ 
\midrule
\endhead
\midrule \multicolumn{6}{r}{\textit{Continued on next page}} \\
\endfoot
\endlastfoot
\textnumero 1.3$\star$	& 1		& 6		     & 	 4		  & 20                              & $V_{2,3}\subset \P^5$ smooth complete intersection of a quadric and a cubic
\\
\textnumero 1.5$\star$	& 1		& 10		     & 	 6		  & 10                              & Gushel–Mukai 3-fold
\\ 
& & & & &  (i) section of Plücker embedding of $\Gr(2,5)$ by codimension 2 subspace and a quadric 
\\ 
& & & & &  (ii) double cover of \textnumero 1.15 with branch locus an anticanonical divisor
\\
\textnumero 1.6 & 1		& 12		     & 	 7		  & 7                              & $X_{12}\subset \P^{8}$, section of $O\Gr_+(5, 10)\subset \P^{15}$ by a linear subspace of codimension 7
\\
\textnumero 1.8	& 1		& 16		     & 	 9		  & 3                             & $X_{16}\subset \P^{10}$, section of $\LGr(2,5)\subset \P^{13}$ by a linear subspace of codimension 3
\\
\textnumero 1.9 & 1		& 18  	     & 	 10		  & 2                              & $X_{18} \subset \P^{11}$, section of $G_2/P\subset \P^{13}$ by a linear subspace of codimension 2
\\
\textnumero 1.10 & 1		& 22		     & 	 12		  & 0                             & $X_{22}$, zero locus of three sections of $\Lambda^2\mathcal{U}^{\vee}$, where $\mathcal{U}$ is the universal sub-bundle on $\Gr(3, 7)$ 
\\
\textnumero 1.14 & 2		& 4		     & 	 		  & 2                            & $V_4\subset \P^5$, smooth complete intersection of two quadrics
\\
\textnumero 1.15	& 2		& 5		     & 	 		  & 0                             & $V_5\subset \P^6$, section of $\Gr(2,5)\subset \P^9$ by a linear subspace of codimension 3
\\
\textnumero 1.16 & 3		& 2		     & 	 		  & 0                             & $Q^3\subset \P^4$, smooth quadric
\\
\textnumero 1.17	& 4		& 1		     & 	 		  & 0                            & $\P^3$
\\[1mm]
\bottomrule
\multicolumn{6}{c}{\textrm{Table~\ltables{tab1}. Minimal geometrically rational real Fano threefold with $\rho(X_\C)=1$}}
\\
\end{longtable}
\end{center}

\proof
Looking at \cite[Section 12.2]{IP99}, we can list the families of rational complex Fano threefolds, which appear in Table~\ref{tab1}. 
\endproof

In general, varieties satisfying \eqref{assumptions} do not verify $\rho(X_\C)=1$ (i.e.\ they are in general not geometrically minimal). Still, Prokhorov classified in \cite{P13} complex Fano threefolds of high Picard rank endowed with an action of a finite group $G$ on the Picard group which preserves the intersection form and the anticanonical class and such that $\rho(X)^G=1$ (see also \cite{CFST16,CFST18}). When $G=\Gal(\C/\R)\cong\Z/2\Z$, we obtain the following result.

\begin{lemma}\label{lemma_tab_2}
Let $X$ be a smooth real Fano threefold verifying \eqref{assumptions} and such that $\rho(X_\C)>1$. Then $X_\C$ belongs to one of the families in Table~\ref{tab2}. 
\end{lemma}

\proof
From the eight families appearing in the table in \cite{P13} we can eliminate (1.2.1) and (1.2.5), which are not rational (see \cite{AB92}). Moreover (1.2.6), (1.2.7) and (1.2.8) are excluded because $G$ has order 2 (see the proof of \cite[Proposition~4.3 - Lemma 4.4]{P13}).
\endproof

\begin{center}
\phantomsection\refstepcounter{tables}
\renewcommand{\arraystretch}{1.2}
\begin{longtable}{@{}l|c|c|c|c|p{9.5cm}@{}}
\toprule
\textnumero $m.n$		& $\iota$	  & $d$ & $g$ & $h^{1,2}$	 &  Description of $X_\C$
\\
\midrule
\endfirsthead
\toprule
\textnumero $m.n$		& $\iota$	  & $d$ & $g$ & $h^{1,2}$	 &  Description of $X_\C$
 \\ 
\midrule
\endhead
\midrule \multicolumn{6}{r}{\textit{Continued on next page}} \\
\endfoot
\endlastfoot
\textnumero 2.12	& 1			     & 	20		  & 11     &     3                  & $X_{(3,3)}$, smooth intersection of three divisors of degree $(1,1)$ in $\P^3\times \P^3$
\\
\textnumero 2.21	& 1			     & 	28 		  & 15       &     0                & $X_{(4,4)}$, blow up of $Q^3\subset \P^4$ along a smooth rational quartic curve
\\
\textnumero 2.32	& 2		     & 	 6	  & 25       &    0                   & $X_{(2,2)}$,  smooth divisor of degree $(1,1)$ in $\P^2\times \P^2$
\\[1mm]
\bottomrule
\multicolumn{6}{c}{\textrm{Table~\ltables{tab2}. Minimal geometrically rational real Fano threefold with $\rho(X_\C)>1$}}
\\
\end{longtable}
\end{center}

We need to recall some results about rationality of geometrically rational Fano threefolds. Many results have been obtained for special families of smooth Fano threefolds over $\R$ and more generally over non-closed fields, via a the study of refined obstruction to rationality introduced by Benoist-Wittenberg and developed by Hassett-Tscinkel and Kuznetsov-Prokhorov (cf.\ \cite{CTSSD87}, \cite{BW23,BW20}, \cite{HT21,HT21b,HT19},\cite{KP23,KP24}).
We restate here the results about the families in Tables \ref{tab1} - \ref{tab2}, when $k=\R$. In the following,
$\Fa_d(X)$ denotes the Hilbert scheme of degree $d$ 
genus zero curves on $X$.

\begin{theorem}[Benoist-Wittenberg, Hassett-Tscinkel, Kuznetsov-Prokhorov]
\label{thm_KP} \

Let $X$ be a smooth real Fano threefold verifying \eqref{assumptions}.

\begin{enumerate}
\item If $X_\C$ belongs to family \textnumero 1.15, then $X$ is $\R$-rational.
\item If $X_\C$ belongs to family \textnumero 1.6, \textnumero 1.10, \textnumero 1.16 or \textnumero 1.17, then 
$$\text{$X$ is $\R$-rational $\iff$ $X(\R)\neq \emptyset$.}$$
\item If $X_\C$ belongs to family \textnumero 1.14, then
$$\text{$X$ is $\R$-rational $\iff$ $X(\R)\neq \emptyset$ and $\Fa_1(X)(\R)\neq \emptyset$}.$$
\item If $X_\C$ belongs to family \textnumero 1.8, then 
$$\text{$X$ is $\R$-rational $\iff$ $X(\R)\neq \emptyset$ and $\Fa_3(X)(\R)\neq \emptyset$}.$$
\item If $X_\C$ belongs to family \textnumero 1.9, then
$$\text{$X$ is $\R$-rational $\iff$ $X(\R)\neq \emptyset$ and $\Fa_2(X)(\R)\neq \emptyset$}.$$
\item If $X_\C$ belongs to family \textnumero 2.21 or \textnumero 2.32, then 
$$\text{$X$ is $\R$-rational $\iff$ $X(\R)\neq \emptyset$.}$$
\item If $X_\C$ belongs to family \textnumero 2.12, then $X$ is never $\R$-rational.
\end{enumerate}
Moreover, for all previous families, $X$ is $\R$-unirational $\iff$ $X(\R)\neq \emptyset$.
\end{theorem}

\proof
The statement is a special case of \cite[Theorem~1.1]{KP23} and \cite[Theorem~1.2]{KP24}.
\endproof

\begin{corollary}
Let $X$ be a smooth real Fano threefold such that $X(\R)\neq\emptyset$ and $X_\C$ belongs to family \textnumero 1.6, \textnumero 1.10, \textnumero 1.15, \textnumero 1.16, \textnumero 1.17, \textnumero 2.21 or \textnumero 2.32. Then $X$ is rational and, as a consequence, $X(\R)$ is connected.
\end{corollary}

\proof
This is a direct consequence of the previous result and the fact that the number of connected components of the real locus is a birational invariant for smooth projective varieties, see e.g. \cite[Theorem~2.3.12]{M20}.
\endproof

An easy observation is the following.

\begin{corollary}
\label{cor:X18}
Let $X$ be a smooth real Fano threefold such that $X(\R)\neq\emptyset$ and $X_\C$ belongs to family \textnumero 1.9. Then $X$ is rational and, as a consequence, $X(\R)$ is connected.
\end{corollary}

\proof
By Theorem~\ref{thm_KP}, it is enough to prove that $X$ contains a real conic, but this is clear, since the general point of $X_{\C}$ has 9 conics passing through it by \cite[Table 2.8.1]{T89}. Since $X$ is unirational, there exists a real conic passing through any general $x\in X(\R)$.
\endproof

All members of families \textnumero 1.3 and 1.5 are conjecturally geometrically irrational, so they play a marginal role in our analysis. Still, we provide few examples.

\begin{proposition}
\label{prop_1.3-5}
There exists smooth real Fano threefolds belonging to family \textnumero 1.3 whose real loci have exactly $s$ connected components, for $s=1,2$ and others belonging to family \textnumero 1.5 whose real loci have exactly $s$ connected components, with $1\le s \le 10$.
\end{proposition}

\proof
Let $Y\subset \P^5$ be a smooth real cubic  with non connected real locus and let $P,P'$ be two real points belonging to distinct connected components of $Y(\R)$. Let $Q\subset \P^5$ be a real quadric with nonempty real locus passing through $P,P'$. We can choose $Q$ such that the intersection $X=Y\cap Q$ is smooth and then $X$ is a real Fano threefold $X$ in family  \textnumero 1.3 whose real locus has two connected components. A similar construction provides an example with one connected component.

Let $\mathrm{Gr}(2,5)\subset\mathbb{P}^9$ be the Grassmannian of $2$-planes in $\C^5$ in its Pl\"ucker embedding, and $V$ be a smooth intersection of $\mathrm{Gr}(2,5)$ with a linear subspace of codimension~$3$. Then 
$V$ is the unique smooth complex Fano threefold in the deformation family \textnumero 1.15. A general member $S$ of $|-K_{V}|$ is a smooth K3 surface. Let $U\to V$ be the double cover branched over $S$. 
Then $U$ is a smooth Fano threefold that belongs to the deformation family \textnumero 1.5, actually one of the \emph{special} members of this family \cite[Theorem~1.1]{De20}. 
Any real forms of $V$ is rational, by Theorem~\ref{thm_KP} (see also \cite{DK19} for the geometry of these real forms). Choose one of them and denote it by $Z$; thus $Z(\R)$ is connected. Let $T\in |-K_{Z}|$ be a real K3 surface such that $T(\R)$ has $10$ connected components. To prove the existence of such a surface in this linear system we can follow the same lines as \cite{Kh76}, see also \cite[VIII.4]{S89}. 

The double cover $W\to Z$ branched over $T$ has two real forms exchanged by the deck involution: one of them whose real locus is connected and the other whose real locus has $10$ connected components. 
Similar constructions gives examples with all intermediate number of connected components.
\endproof

In the rest of this section, we will study connectedness of $X(\R)$ for the remaining families, namely \textnumero 1.14, \textnumero 1.8,  \textnumero 2.12. 

\subsection{Family \textnumero 1.14: Krasnov's classification}
\label{sec_V4}

We recall here some of the results from \cite{K18}, where the author classifies topological types for the real locus of three-dimensional smooth complete intersection of two quadrics over $\mathbb{R}$. 

Unirationality is well known (see \cite[Proposition~2.3]{CTSSD87}), while the problem of rationality has been studied in \cite{HT21b} and completely solved in \cite[Theorem~A]{BW23} (see the part of Theorem~\ref{thm_KP} on $V_4$).

In order to understand the results in \cite{K18}, we need to recall here the notion of {\it isotopy classes}. Here we follow \cite[Section~11]{HT21b}.

Let $X \subset \mathbb{P}^5$ smooth complete intersection of two quadrics over $\mathbb{R}$, defined by $q_0$ and $q_1$. One can associate a pencil $Y\subset \mathbb{P}^5\times \mathbb{P}^1$ defined by 
$$\lambda_0q_0 + \lambda_1q_1 = 0.$$ 
Take the $2:1$ covering $\gamma\colon S^1 \to \mathbb{P}^1$ and the base-change 
\[\xymatrix {
   \tilde{Y} \ar[r] \ar[d] &  Y \ar[d] \\
  S^1 \ar[r] & \mathbb{P}^1 \\
}\]
Now $\tilde{Y}$ defines a well-defined family of quadratic forms over $S^1$ and one can take the positive index of inertia 
$$I^+\colon S^1\to \mathbb{Z}$$
computing the number of positive eigenvalues. By construction $I^+$ is piecewise constant and jumps $2k\le 12$ times, with height $\pm 1$ (here we use that $X$ is smooth). 
\\One says that a point of discontinuity for $I^+$ is {\it positive} if the value of the inertia increases by one as we cross it anti-clockwise. This produces a partition 
\begin{equation}\label{odd_dec}
k=k_1+k_2+\cdots+k_{2s+1}
\end{equation}
given by the numbers of consecutive positive points of discontinuity, moving anticlockwise on $S^1$. One also quotients out the set of these decompositions by cyclic permutations or reversal of the order of the sum.

\begin{definition}
Let $\mathcal{B}_3$ denote the space of three-dimensional smooth complete intersection of two quadrics over $\mathbb{R}$. Then the connected components of $\mathcal{B}_3(\mathbb{R})$ are called the {\it (rigid) isotopy classes}.
\end{definition}

\begin{proposition}[Krasnov]
Isotopy classes of smooth three-dimensional complete intersections of two quadrics correspond to equivalence classes of odd decompositions \eqref{odd_dec}, where $0 \le k\le 6$ is even. 
\end{proposition}

So we will refer to an isotopy class with its corresponding partition \eqref{odd_dec}: either $(0)$ or $(k_1,\ldots,k_{2s+1})$ with $0<k_1\le k_2\le\cdots \le k_{2s+1}$.

The following result is a direct consequence of the results in \cite{K18}.

\begin{proposition} [Krasnov]
\label{proposition:Kra}
Let $X$ be a real form of $V_{4}$, then 
$$X(\R) \  \text{is disconnected } \iff \  \#\pi_0(X(\R))=2 \iff \ X \text{ has isotopy class (1,1,4)}.$$
\end{proposition}

\proof
If $X$ has isotopy class $(0)$, then $X(\R)=\emptyset$. The classification result in \cite[Theorem~5.4]{K18} implies the result.
\endproof


\subsection{Family \textnumero 1.8: weak Fano model with quadric fibration}
\label{sec_X16}

In this part we produce examples in family \textnumero 1.8 whose real locus has several connected components, building on \cite[(2.3.8)]{T22} and \cite[Example 4.9]{Chelbook23}. We state here a key lemma to control the number of connected components of the real locus under smoothing.  We follow \cite{Na97}.

\begin{lemma}
\label{lemma:smoothing}
    Let $X_0$ be a real singular Fano threefold whose singularities are ordinary double points. 
   Then $X_0$ is smoothable over $\mathbb{R}$. In particular, there exists a real smooth Fano threefold $X$ of the same degree and same complex Picard number whose real locus $X(\mathbb{R})$ has at least the same number of connected components as $X_0(\R)$.
 Furthermore, if the real locus is isomorphic locally around all ordinary double points to the cone over $S^1\times S^1$, then $\#\pi_0(X(\R))=\#\pi_0(X_0(\R))$.
\end{lemma}
\proof
After base-change, we have $H^2(X_{0,\C},T_{X_{0,\C}})=0$ 
by \cite[Proposition~4]{Na97}, and 
\cite[Theorem~11]{Na97} implies that $X_{0,\C}$ is smoothable by a flat deformation and  by \cite[Proposition 3 and Lemma 12]{Na97}, the Kuranishi space $\operatorname{Def}(X)$ is smooth and universal.  
Since $X_0$ is defined over $\mathbb{R}$, $\operatorname{Def}(X)$ is endowed with an anti-linear involution whose fixed locus parametrizes real infinitesimal deformations of $X_0$. Hence $X_0$ is smoothable over $\mathbb{R}$ and we get a real smooth Fano threefold $X$ of the same degree. By Jahnke-Radloff \cite[Theorem~1.4]{JR11}, the Picard number of $X_\C$ is the same as the Picard number of $X_{0,\C}$.

The cone over $S^2$ admits two different smoothing, one of which locally disconnects the real locus. On the other hand, the cone over $S^1\times S^1$ admits only one smoothing which is connected.
\endproof

\begin{proposition}
\label{prop:X16}
There exists a smooth real Fano threefold $X$ belonging to family \textnumero 1.8 whose real locus has $3$ connected components.
\end{proposition}

\proof 
Let $W =\mathbb{P}(\mathcal{O}_{\mathbb{P}^1} \oplus\mathcal{O}_{\mathbb{P}^1} \oplus\mathcal{O}_{\mathbb{P}^1}(1) \oplus\mathcal{O}_{\mathbb{P}^1}(1) )$ and let $\pi\colon W\to \mathbb{P}^1$ be the natural projection. Denote by $H$ the tautological bundle and by $F$ a fibre of $\pi$. Write $t_0,t_1$  for the coordinates on $\mathbb{P}^1$ and $x,y,z,w$ be coordinates on the fibre with $x,y$ sections of $H$ and $z,w$ sections of $H-F$.

Let $V$ be the divisor in $\vert 2H+F\vert$ given by the following equation
$$
(t_0+t_1)x^2+(t_0+2t_1)y^2+t_1(t_1-t_0)(t_1-2t_0)z^2+t_0(t_1-3t_0)(t_1-4t_0)w^2=0
$$

The manifold $V$ is a Picard-rank-two weak Fano threefold with a fibration $\pi:=\pi_{\vert V}\colon V \to \mathbb{P}^1$. We have $(-K_V)^3=16$ and the real locus $V(\mathbb{R})$ has $3$ connected components over the intervals $t_0=1$, $1\leq t_1\leq 2$, $3\leq t_1\leq 4$ and $\infty\leq t_1\leq 0$ of $\mathbb{P}^1(\mathbb{R})$. Here $t_1=\infty$ is for the point $[t_0,t_1]=[0:1]$.
Let $C$ be the curve on $V$ given by $\{z=w=0\}$. 
The anticanonical map of $V$ is small, and $C$ is the only curve with trivial intersection $-K_V\cdot C = 0$.
The curve $C$ is smooth rational and is a bisection of $\pi$. The curve $C$ is the complete intersection of the two surfaces $S_1:\{z=0\}$ and $S_2:\{w=0\}$ in $V$. It is smooth along each plane $S_i$ then it is a $(-1)$-curve in each plane $S_i$, thus its normal bundle $\mathcal{N}_{C/V}$ is isomorphic to $\mathcal{O}_{\mathbb{P}^1}(-1) \oplus\mathcal{O}_{\mathbb{P}^1}(-1)$.

The real locus $C(\mathbb{R})$ meets only one connected component of $V(\mathbb{R})$, namely the one over the interval $t_0=1,-2\leq t_1\leq -1$ of $\mathbb{P}^1(\mathbb{R})$.

Let $X_0$ be the singular Fano threefold obtained by contracting $C$, then $X_0$ is a singular Fano threefold of index $1$ and genus $9$ whose unique singularity is an ordinary double point whose real locus is the cone over $S^1\times S^1$. Its real locus $X_0(\R)$ has the same number of connected components as $V(\R)$.
By Lemma~\ref{lemma:smoothing}, we get a real smooth Fano threefold $X$ in family \textnumero 1.8 whose real locus $X(\mathbb{R})$ has exactly $3$ connected components. 
\endproof

\begin{remark}
With the same construction, we can produce further examples with two connected components starting with $V$ given by the equation
$$
(t_0+t_1)x^2+(t_0+2t_1)y^2+t_1(t_1-t_0)(t_1-2t_0)z^2+t_0(t_1^2+t_0^2)w^2=0
$$
or one connected components with the equation
$$
(t_0+t_1)x^2+(t_0+2t_1)y^2+t_1(t_1^2+2t_0^2)z^2+t_0(t_1^2+t_0^2)w^2=0.
$$
\end{remark}


\subsection{Family \textnumero 2.12: minimal irrational members with connected and disconnected real loci.}
\label{sec_X33}

In this part we produce examples of {\it minimal} smooth real Fano threefold belonging to family \textnumero 2.12, whose real locus is connected or disconnected. We recall that Theorem~\ref{thm_KP} implies that any such member is $\R$-irrational. 

We recall that over $\C$ the elements of family \textnumero 2.12 can be ralised as smooth intersection of three divisors of degree $(1,1)$ in $\P^3\times \P^3$, so, in order to construct explicit examples, we need to recall some well-known facts on real forms of $\P^3\times\P^3$.

Any real structure on $\P^3$ is equivalent to one of the following.
\begin{enumerate}
\item 
the complex conjugation
$$
\sigma_0\colon [x_0:x_1:x_2:x_3]\mapsto [\bar x_0:\bar x_1:\bar x_2:\bar x_3]
$$
whose associated real locus is $\P^3(\R)$
\item the one given by 
$$
\sigma_1\colon [x_0:x_1:x_2:x_3]\mapsto [-\bar x_1:\bar x_0:-\bar x_3:\bar x_2]
$$
whose associated real locus is empty.
\end{enumerate}

We denote by $\sigma_{twist}$ the real structure on $\P^3\times\P^3$ given by
$$
([x_0:x_1:x_2:x_3],[y_0:y_1:y_2:y_3])\mapsto ([\bar y_0:\bar y_1:\bar y_2:\bar y_3],[\bar x_0:\bar x_1:\bar x_2:\bar x_3])\;.
$$

\begin{lemma}
\label{lemma:sigma_twist}
Any real structure on $\P^3\times\P^3$ is equivalent to one of the following

\begin{enumerate}
\item 
$
\sigma_0\times \sigma_0
$
whose associated real locus is $\P^3(\R)\times\P^3(\R)$.
\item 
$
\sigma_0\times\sigma_1
$
whose associated real locus is empty.
\item 
$
\sigma_1\times\sigma_1
$
whose associated real locus is empty.
\item
$
\sigma_{twist}
$
whose real locus is diffeomorphic to the  $6$-dimensional real manifold underlying $\P^3(\C)$. 
\end{enumerate}
\end{lemma}

\proof
The four real structures above are pairwise nonequivalent. Conversely, we have $\operatorname{Aut}_\C(\P^3\times\P^3)\simeq \operatorname{PGL}_4(\C)^2\rtimes \Z/2$ and the Galois cohomology pointed set $H^1(G,\operatorname{Aut}_\C(\P^3\times\P^3))$ 
has four elements,  see e.g. \cite{GS17} for details.
\endproof

\begin{remark}
The  $6$-dimensional real manifolds $\P^3(\R)\times\P^3(\R)$ and $\P^3(\C)$ are orientable.
\end{remark}

Starting form \cite[Example 5.4]{Chelbook23}, we show that this construction provides an interesting example whose real locus has one connected component.

\begin{proposition}
\label{prop:X33-connected}
There exists a minimal (hence $\R$-irrational) smooth real Fano threefold in family \textnumero 2.12 with connected real locus.
\end{proposition}

\proof
Let $\mathcal{C}\subset \P^3$ be the unique $\operatorname{PSL_2(\mathbf{F}_7)}$-invariant smooth curve of degree $6$ and genus $3$ which is real by unicity.
Let $\pi\colon X \to \P^3$ be the blow-up of $\mathcal{C}$.
Then the (non-minimal) threefold $X$ belongs to family \textnumero 2.12, and can be described in $\P^3\times \P^3$ by (see \cite[(5.4.3)]{Chelbook23})
$$
\begin{cases}
 x_0y_1+x_1y_0-\sqrt2x_2y_2=0\\
x_0y_2+x_2y_0-\sqrt2x_3y_3=0\\
x_0y_3+x_3y_0-\sqrt2x_1y_1=0
\end{cases}
$$

The involution $\tau\in \operatorname{Aut}(\P^3\times \P^3)$ given by 
$$
([x_0:x_1:x_2:x_3],[y_0:y_1:y_2:y_3])\mapsto ([y_0:y_1:y_2:y_3],[x_0:x_1:x_2:x_3])
$$
leaves $X$ invariant hence induces an element of $\operatorname{Aut}(X)$ which we still denote by $\tau$ (this involution is denoted $\sigma$ in \cite[Example 5.4]{Chelbook23}).

Then the restriction to $X_\C$ of the real structure $\sigma_{twist}$ on $\P^3\times\P^3$ defined before Lemma~\ref{lemma:sigma_twist} is also the one obtained by composition of $\tau$ and the canonical real structure $\sigma_0$ corresponding to the complex conjugation $\sigma_0$ on $\P^3$, lifted through $\pi$. We denote by $\sigma$ this real structure and by $Y$ the corresponding real form of $X$. 

Let $H$ be an hyperplane in $\P^3$. By construction, the real Picard number of $X$ is $2$ and $\Pic(X_\C)^{\sigma_0}$ is generated by $\pi^*H$ and the $\pi$-exceptional surface $E$.
The involution $\sigma$ acts on $\Pic(X_\C)$ as follows:
$$
\begin{cases}
\sigma^*(E)\sim 8\pi^*(H)-3E\;,\\
\sigma^*(\pi^*(H))\sim 3\pi^*(H)-E\;.
\end{cases}
$$

Hence $-K_X=\pi^*(H)+\sigma^*(\pi^*(H))$ and $\Pic^\sigma(X_\C) = \Z[-K_X]$
thus the real Picard number of $Y$ is $1$.

The real locus $Y(\R)$, is given by the following equations
$$
\begin{cases}
y_1=\bar x_1\\
y_2=\bar x_2\\
y_3=\bar x_3\\
2\Re(x_0\bar x_1)-\sqrt 2 \abs{x_2}^2=0 \\
2\Re(x_0\bar x_2)-\sqrt 2 \abs{x_3}^2=0 \\
2\Re(x_0\bar x_3)-\sqrt 2 \abs{x_1}^2=0 
\end{cases}
$$
Those ones provide $9$ equations in real variables.

We show now that $Y(\R)$ is connected. First we observe that $Y(\R)$  is contained in the affine chart $\{x_0\ne 0\}\times\{y_0\ne 0\}$. Indeed, if $([x_0:x_1:x_2:x_3],[y_0:y_1:y_2:y_3])$ belongs to $Y(\R)$, and if $x_0=0$, then $x_j=0$ for $j=1,2,3$. Thus we can assume $x_0=1=y_0$. 
The first $6$ real equations define a real affine subspace $V\subset \A^3\times \A^3$ of real dimension $6$.
Letting $x_j=a_j+ib_j$, $a_j,b_j\in\R$ for $j=1,2,3$, we get that $Y(\R) \subset V$ is defined by the equations
$$
\begin{cases}
2a_1-\sqrt 2 (a_2^2+b_2^2)=0 \\
2a_2-\sqrt 2 (a_3^2+b_3^2)=0 \\
2a_3-\sqrt 2 (a_1^2+b_1^2)=0 
\end{cases}\;.
$$

Let $P=(a_1,b_1,a_2,b_2,a_3,b_3)\in Y(\R)$ be a real point and for $t\in [0,1]$ define $P(t)=(a_1(t),b_1(t),a_2(t),b_2(t),a_3(t),b_3(t))$ where
$$
\begin{cases}
a_j(t)=(1-t)a_j\quad j=1,2,3\\
b_2(t)=\sqrt{1-t}\sqrt{2a_1-(1-t)\sqrt2 a_2^2}\\
b_3(t)=\sqrt{1-t}\sqrt{2a_2-(1-t)\sqrt2 a_3^2}\\
b_1(t)=\sqrt{1-t}\sqrt{2a_3-(1-t)\sqrt2 a_1^2}
\end{cases}\;.
$$

Then $[0,1]\to Y(\R), t\mapsto P(t)$ is a continuous path  from $P(0)=P$ to $P(1)=(0,0,0,0,0,0)$ such that $P(t)\in Y(\R)$  $\forall t\in [0,1]$. 

Let's prove that $b_2(t)$ is well defined.
Indeed, $P\in Y(\R)$ 
then $0\leq \sqrt 2b_2^2=2a_1-\sqrt 2 a_2^2$ and for $t\in [0,1]$, we get $2a_1-(1-t)\sqrt2 a_2^2\geq 0$. The same proof shows that $b_1$ and $b_3$ are well-defined. We have proven that $Y(\R)$ is path connected.
\endproof

The final proposition of this part provides an example with two components, obtained in the spirit of \cite{CTZ24-pre}.

\begin{proposition}\label{prop:X33}
There exists a smooth real Fano threefold belonging to family \textnumero 2.12 with disconnected real locus.
\end{proposition}

\proof
Let  $Y\subset\P^4$ be a real cubic hypersurface whose singular locus is a pair of conjugated ordinary double points $\{P,\bar P\}\subset Y(\C)\setminus Y(\R)$ and real locus $Y(\R)$ has two connected components.
An explicit example of such a cubic is the following:
$$x_0^2x_2 + x_1^2x_2 + x_2^3 -4x_2^2x_3 + x_0x_4^2 + x_2x_3^2 + x_3^3 + x_2^2x_4 +
 x_2x_3x_4 + x_3^2x_4 + x_1x_4^2 + 7x_2x_4^2 + x_3x_4^2 + x_4^3=0$$ 
 where $[x_0:x_1:x_2:x_3:x_4]$ are coordinates in $\P^4$.

The image of the projection onto $\P^2$ with coordinates $[x_2:x_3:x_4]$ is given by

$-x_2^4 + 4x_2^3x_3 -x_2^2x_3^2 -x_2x_3^3 + \frac14x_3^4 -x_2^3x_4 -x_2^2x_3x_4 -x_2x_3^2x_4 -7 x_2^2x_4^2 -x_2x_3x_4^2 -x_2x_4^3 +
 \frac14x_4^4=0.$
 
 \begin{figure}[ht]
\centering
\includegraphics[height =4cm]{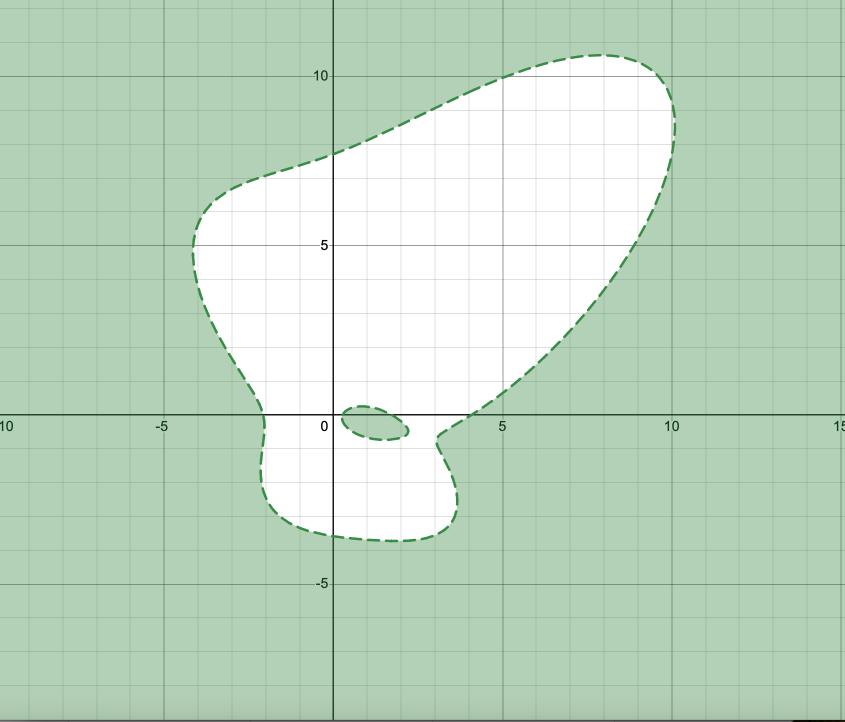}
\caption{Image of the projection of $X(\R)$ on $\p^2$ (in green).}
        \label{fig.cubic}
\end{figure}

 Let $\pi, \bar\pi\colon Y\dasharrow \p^3$ be the conjugated projections of $Y$ from the singular points $P,\bar P$. The product map gives a real rational map $f\colon Y\dasharrow  W$  where $W$ is the real form of $\p^3\times \p^3$ corresponding to the real structure $\sigma_{twist}$ (see Lemma~\ref{lemma:sigma_twist}). 
 The image $X_0=\overline{f(Y)}$ has a unique ordinary double point which is the result of the contraction of the line passing through $P$ and $\bar P$. The real locus $X_0(\R)$ still has two connected components because the real locus of the contracted line is connected. By Lemma~\ref{lemma:smoothing}, we get a real smoothing $X$ of $X_0$ which is a real Fano threefold belonging to family \textnumero 2.12 with disconnected real locus. 
\endproof


\section{Maximal number of connected components and a rationality criterion}
\label{sec:criterion}

In this section we discuss an invariant for real varieties, which encodes the maximal number of connected components in families and prove the main result of this paper. 

\subsection{The invariant $s_X$}
We recall that two smooth complex varieties $Y$ and $Z$ are {\it deformation equivalent}, denoted $Y\overset{\defo}{\sim} Z$ if there exists a finite chain $Y=Y_0,Y_1,\dots,Y_r=Z$ 
 of smooth complex varieties such that $Y_i$ and $Y_{i-1}$ are complex isomorphic to fibers of a smooth morphism over a smooth curve. 
 Note that if $Y\overset{\defo}{\sim} Z$, then $Y(\C)$ and $Z(\C)$ are $\mathcal{C}^\infty$-diffeomorphic, by Ehresmann's theorem.

\begin{definition}\label{def:sx}
Let  $X$ be a smooth real projective variety. Then one defines
$$
s_X:=\max\{\#\pi_0(X'(\R))\,| \,  \,X'_\C\overset{\defo}{\sim} X_\C\}.$$
\end{definition}

Note that $s_X$ is well-defined for any real projective variety by Smith-Thom inequality~\eqref{eq.smiththom}. Note that $s_X$ is in general hard to determine and known only in few cases. Here are some examples.

\begin{examples} 
\label{example:sx}
\begin{enumerate}
\item Let $X,X'$ be real smooth projective curves, then $X'_\C\overset{\defo}{\sim} X_\C$ if and only if they have the same genus $g$. Moreover $s_X=g+1$ by Harnack's Theorem.

\item 
\label{item:dp}
Let $S_d$ be a del Pezzo surface of degree $d$, then $s_{S_1}=5$, $s_{S_2}=4$, $s_{S_3}=s_{S_4}=2$, for $5\leq d \leq 9$, $s_{S_d}=1$, see \cite[Chapter 4]{M20}.
    \item 
Let $X\subset\P^N$ and $X'\subset\P^N$ be real smooth hypersurfaces of the same degree, then $X'_\C\overset{\defo}{\sim} X_\C$. 
\item Let $X_d\subset \mathbb{P}^3$ be a degree-$d$ real smooth surface, then $s_{X_1}=s_{X_2}=1$, $s_{X_3}=2$, $s_{X_4}=10$, but $s_{X_5}$ is not known.
\item Let $N$ be an odd integer, and $X=\P^N$. There exists a unique non-trivial Severi-Brauer variety $X'$, and by definition $X'_\C\overset{\defo}{\sim} X_\C$.
\item Let $X=\P^1\times \P^1$ and $X'=\P^1\times C$ be the product of $\P^1$ with a pointless conic $C$, then $X'$ does not admit a real quadric model in $\P^3$ but $X'_\C\overset{\defo}{\sim} X_\C$.
\item Let $X$, $X'$ be real smooth Fano threefolds such that $X_\C$ and $X'_\C$ belong to the same family then $X'_\C\overset{\defo}{\sim} X_\C$.
\end{enumerate}
\end{examples}

We stress that the number $s_X$ is not a birational invariant (see Remark~\ref{bir_inv}).
Moreover, in what follows, we will always assume that $\xr\ne\varnothing$. For the classification of families containing a member $X$ such that $X(\R)=\varnothing$, see \cite[Proposition~B]{ackm3}.

\subsection{A rationality criterion}

We state now the main result of this paper.

\begin{theorem}
\label{thm_newmain}
Let $X$ be a smooth geometrically rational real Fano threefold  with $\xr\ne\varnothing$.
If $s_X=1$, 
then $X$ is rational.
\end{theorem}

We introduce the following notation.
\begin{definition}
\label{dfn.smn}
Let $X$ be a geometrically rational smooth real Fano threefold  belonging to family \textnumero $m.n$. Then one denotes $s_{m.n}:=s_X.$
\end{definition}

Our strategy consists in going through all families of (sometimes conjecturally) rational complex Fano threefolds and proving that either $s_{m.n}>1$ or the implication $[\xr\ne\varnothing\implies X\ \text{is rational}]$ holds for every real form of a variety belonging to family \textnumero $m.n$.

\begin{lemma}
\label{lem_dp}
Let $X$ be a smooth real Fano threefold such that
$X_\mathbb{C}\cong\mathbb{P}^1\times S_d$ where $S_d$ is a smooth del Pezzo surface of degree $d$ $(S_8=\F_1$ or $S_8=\P^1\times\P^1)$. Then
\begin{enumerate}
\item 
$X$ is rational if and only if its real locus is nonempty and connected. 
\item
 Let $s_d$ be the maximal number of connected components of the real locus of a smooth real del Pezzo surface of degree $d$. If $X$ belongs to family \textnumero $m.n$, then $s_{m.n}=s_{11-m}$.
\end{enumerate}
\end{lemma}

\proof
As $X_\C$ is the product of a curve and a surface, $X$ is as well a product, hence $X\cong U\times V$ where $U$ is a real form of $\P^1$ and $V$ is a real form of a del Pezzo surface 
of degree $10-m$ where $m=\rho(X_\C)$. Then $X$ is rational if and only if $U=\P^1$ and $V$ is rational. By Comessatti's Theorem,
we get that $X$  is rational if and only if $V(\R)$ is connected if and only if $X(\R)$ is connected.
\endproof

\subsection{The proof of Theorem \ref{thm_newmain}}

In what follows, let $X$ be a geometrically rational smooth real Fano threefold. We will study families depending on the Picard rank.

\subsubsection{$\rho(X_\C)=1$}

\begin{proposition}
\label{prop_m=1}
Let $X$ be a geometrically rational smooth real Fano threefold belonging to family \textnumero $1.n$. Assume that $\xr\ne\varnothing$. Then $X$ is rational and $s_{1.n}=1$, unless $n\in\{3,5,8,14\}$.
\end{proposition}
\proof
There are 17 families of smooth Fano threefolds with Picard number $1$.

For $n\in\{1,2,4,7,11,12,13\}$, no $X$ belonging to family \textnumero $1.n$ are geometrically rational.

For $n\in\{3,5\}$, the general member in family \textnumero $1.n$ is not geometrically rational and all members are conjecturally not geometrically rational. Any $X$ belonging to family \textnumero $1.n$, with $n\in\{6,8,9,10,14,15,16,17\}$,  is geometrically rational.

Let $X$ belonging to family \textnumero $1.n$ such that $\xr\ne\varnothing$. By Theorem~\ref{thm_KP} 
for $n\in\{6,10,15,16,17\}$ and by Corollary~\ref{cor:X18} for $n=9$, $X$ is rational.
\endproof

\begin{proposition}
\label{rang_1_more}
The following holds:
\begin{itemize}
\item $s_{1.3}\geqslant 2$;
\item $s_{1.5}\geqslant 10$;
\item $s_{1.8}\geqslant 3$ ;
\item $s_{1.14}=2$.
\end{itemize}
\end{proposition}

\proof
This is a consequence of Propositions~\ref{prop_1.3-5}, \ref{prop:X16} and \ref{proposition:Kra}.
\endproof

\subsubsection{$\rho(X_\C)=2$}

\begin{proposition}
\label{prop_m=2}
Let $X$ be a geometrically rational smooth real Fano threefold belonging to family \textnumero $2.n$. Assume that $\xr\ne\varnothing$. Then $X$ is rational and $s_{2.n}=1$, unless $n\in\{10,12,16,18\}$.
\end{proposition}
\proof
There are 36 families of smooth Fano threefolds with Picard number $2$.

For $n\in\{1,2,3,5,6,8,11\}$, no $X$ belonging to family \textnumero $2.n$ are geometrically rational. Any $X$ belonging to other family \textnumero $2.n$ is geometrically rational.
If $n\in\{4,7,9,13, 14,15,17,19,20,22,23,25,26,27,28,29,30,31,33,35\}$, $X_\C$ admits an extremal birational contraction $f\colon X_\C \to Y$, where $Y$ is $\P^3_\C$, the quadric $Q\subset \P^4_\C$, or the quintic $V_5\subset \P^6_\C$ (see \cite{MM86,MM03,Mat-erratum}). By \cite[Lemma~2.5]{ACKM24}, there exists a birational morphism $g\colon X \to W$, where $W$ is a real form of $Y$. Since $\xr\ne\varnothing$, we deduce that $W(\R)\ne\varnothing$, and therefore $W$ is rational over $\R$. We deduce the statement for these values of $n.$ 
\\Assume $n=21$: if $\rho(X)=1$, then $X$ is rational by Theorem~\ref{thm_KP}, if $\rho(X)=2$, then $X$ admits an extremal birational contraction onto a quadric $Q\subset \P^4_\C$ with $Q(\R)\ne\varnothing$. We conclude that $s_{2.21}=1$.
\\Assume $n\in\{24,34,36\}$: $X_\C$ has an extremal contraction that produces a $\P^1$-bundle $f\colon X_\C \to \P^2_\C$, while the second extremal ray corresponds to a conic bundle for $n=24$, a $\P^2$-bundle for $n=34$ and a divisorial contraction for $n=36$, respectively. This implies that the Galois action cannot exchange the two rays and, as a consequence, that $f$ descends to a $\P^1$-bundle structure over $\R$. We conclude that $X$ is rational and, therefore, $s_{2.n}=1$.
\\Assume $n=32$: if $\rho(X)=1$, then $X$ is rational by Theorem~\ref{thm_KP}, otherwise, it admits two $\P^1$-bundle structure over $\P^2$, hence $X$ is rational and $s_{2.32}=1$.
\endproof

We study now the remaining families with geometric Picard rank two.

\begin{proposition}\label{rang_2_more}
The following holds:
\begin{itemize}
\item[(i)] $s_{2.10}=2$;
\item[(ii)] $s_{2.12}\geq 2$;
\item[(iii)]$s_{2.16}\ge 2$;
\item[(iv)] $s_{2.18}= 3$.
\end{itemize}
\end{proposition}

\proof
Let $X$ be a geometrically rational smooth real Fano threefold belonging to family \textnumero 2.10. By \cite{MM86,MM03,Mat-erratum} and \cite[Lemma~2.5]{ACKM24}, there exists a birational morphism $f\colon X \to Y$, where $Y$ is a real form of a quartic $V_4$ and  $f$ is the blow-up of a genus-one curve, intersection of two hyperplanes. This implies that $s_{2.10}\le 2$. To show the equality, let $Y'$ be a real form of a quartic $V_4$ such that $Y'(\R)$ has two connected components. Now take two general rational points on the two connected components and intersect $Y'$ with two hyperplanes passing through those rational points: this produces a genus-one curve $C'$ in $Y'$ and let $X'$ be the blow up of $Y'$ along $C'$. By construction, $X'(\R)$ has two connected components. This shows that $s_{2.10}=2$.
\\The estimate on $s_{2.12}$ is a consequence of Proposition~\ref{prop:X33}. 
\\We produce now an example in family \textnumero 2.16 with two connected components. Let $Q_1, Q_2\subset\P^5$ be the quadric threefolds given by the equations 
$$x_1^2+x_2^2+x_3^2+x_4^2+x_5^2-x_0^2=0$$
$$\alpha^2(x_1^2+x_2^2+x_3^2+x_4^2)+\beta^2x_5^2-x_0^2=0$$ 
where $1<\alpha$ and $0<\beta<1$. Let $Y=Q_1\cap Q_2$ be their intersection. For general $\alpha$ and $\beta$ the threefold $Y$ is smooth. By construction, $Y$ belongs to family \textnumero 1.14. Moreover, all real points are contained in the affine chart $x_0\neq 0$, where the equations become
\begin{gather*}
x_5^2=1-(x_1^2+x_2^2+x_3^2+x_4^2)\\
\beta^2x_5^2= 1-\alpha^2(x_1^2+x_2^2+x_3^2+x_4^2)
\end{gather*}
One checks that the hyperplane section $\{x_5=0\}\cap Y$ has no real point, since $1<\alpha$, while each half-space $\{\pm x_5>0\}\cap Y$ contains real points. Hence $Y(\R)$ has two connected components. Moreover, the following holds:
\begin{enumerate}
\item $Y$ is contained in the quadric given by $x_5^2=\frac{\alpha^2-1}{\alpha^2-\beta^2}x_0^2$;
\item Let $P$ be the plane given by $\left\{x_1=x_2=0,\ x_5=\sqrt{\frac{\alpha^2-1}{\alpha^2-\beta^2}}x_0\right\}$. 
The curve $Y\cap P$ is a smooth real conic $C$.
\end{enumerate}
Let $X$ be the blow up of $Y$ along $C$, then $X$ is a real smooth Fano threefold belonging to family \textnumero 2.16 and $X(\R)$ has two connected components. This implies $s_{2.16}\ge 2$.
\\The family \textnumero 2.18 was extensively studied in \cite{FJSVV24} and \cite{JJ23}. The construction given in \cite[Example 4.6]{JJ23} implies that $s_{2.18}\ge 3$. One can adapt the argument in \cite[Lemma 2.7]{JJ23} to deduce that, even when the discriminant of the conic bundle is not smooth, the bound 
 $s_{2.18}\le 3$ holds true. We conclude that $s_{2.18}= 3$.
\endproof

\begin{remark}\label{bir_inv}
The previous proposition shows that $s_X$ is not a birational invariant. If $X$ belongs to family \textnumero 1.14, then $s_X=2$. If $Y$ is the blowup of $X$ in an elliptic curve which is an intersection of 2 hyperplanes, i.e. $Y$ belongs to familly \textnumero 2.10, then $s_Y=2$. 
But if $Y'$ is the blowup of $X$ in a line, i.e. $Y'$ belongs to familly \textnumero 2.19, then $s_{Y'}=1$.
Indeed $X$ is rational if and only if it contains a real line by Theorem~\ref{thm_KP}. 
\end{remark}

\subsubsection{$\rho(X_\C)=3$}

\begin{proposition}
\label{prop_m=3}
Let $X$ be a geometrically rational smooth real Fano threefold belonging to family \textnumero $3.n$. Assume that $\xr\ne\varnothing$. Then $X$ is rational and $s_{3.n}=1$, unless $n\in\{2,3,4\}$.
\end{proposition}
\proof
There are 31 families of smooth Fano threefolds with Picard number $3$.
No $X$ belonging to family \textnumero 3.1 are geometrically rational. Any $X$ belonging to other family \textnumero $3.n$ is geometrically rational.
The strategy consists in checking for all families the type of extremal contractions of $X_\C$, studied in \cite{MM86,MM03,Mat-erratum}. See also \cite{ACKM24} for details on the geometry of the extremal contractions.
\begin{itemize}
\item[{Family \textnumero 3.5}.] $X_\C$ admits exactly three extremal divisorial contractions $f_i\colon X_\C \to Y_i$, where $Y_1\simeq\P^1\times \P^2$, while $Y_2$ and $Y_3$ are not Fano, see \cite[\S~III.3, p.~74]{Mat95}. The Galois action fixes $f_1$, which descends to $g_1\colon X \to W_1$. Since $\xr\ne\varnothing$, we deduce that $W(\R)\ne\varnothing$, and therefore $W$ is rational over $\R$. This implies that $X$ is rational.
\item[{Family \textnumero 3.6}.] $X_\C$ admits exactly two extremal divisorial contractions $f_i\colon X_\C \to Y_i$, where $Y_1$ belongs to family 2.25 and $Y_2$ belongs to family 2.33, see \cite[\S~III.3, p.~75]{Mat95}. Since they cannot be exchanged by the Galois action, they descend to $g_i\colon X \to W_i$, where $W_i$, $i=1,2$, are rational by Proposition~\ref{prop_m=2}. This implies that $X$ is rational.
\item[{Family \textnumero 3.7}.] $X_\C$ admits exactly three extremal divisorial contractions $f_i\colon X_\C \to Y_i$, where $Y_1\simeq Y_2 \simeq \P^1\times \P^2$ and $Y_3$ belongs to family 2.32, see \cite[\S~III.3, p.~76]{Mat95}. The Galois action fixes $f_3$, which descends to $g_3\colon X \to W_3$, where $W_3$ is rational by Proposition~\ref{prop_m=2}. This implies that $X$ is rational.
\item[{Family \textnumero 3.8}.] $X_\C$ admits exactly two extremal divisorial contractions $f_i\colon X_\C \to Y_i$, where $Y_1\simeq \P^1\times \P^2$ and $Y_2$ belongs to family 2.24, see \cite[\S~III.3, p.~77]{Mat95}. Since they cannot be exchanged by the Galois action, they descend to $g_i\colon X \to W_i$, where $W_i$, $i=1,2$, are rational by Proposition~\ref{prop_m=2}. This implies that $X$ is rational.
\item[{Family \textnumero 3.9}.] $X_\C$ is isomorphic to the blow-up of the cone over the Veronese of $\P^2$ in $\P^5$ with center the disjoint union of the vertex and a quartic curve on $\P^2$ and admits exactly two extremal divisorial contractions of the form $f_i\colon X_\C \to Y$, $i=1,2$, where $Y$ is the blow-up of the Veronese cone in a quartic curve, see \cite[\S~III.3, p.~79]{Mat95}. If those two contractions are fixed by the Galois action, they descend to $\R$ and induce birational morphisms $h_i\colon X \to W$, where $W$ is the Veronese cone over $\R$ (here the relative Picard rank for $h_1$ is 2). This implies that $X$ is rational. 
\\If these two extremal contractions $f_1$ and $f_2$ are exchanged by the Galois action. Then we still obtain a birational morphism $g\colon X \to W$ over $\R$, where $W$ is the Veronese cone over $\R$ (in this second case, the relative Picard rank for $g$ is 1). As before, we conclude that $X$ is rational.
\item[{Family \textnumero 3.10}.] $X_\C$ is isomorphic to the blow-up of the quadric $Q\subset \P^4_\C$ in two disjoint conics and admits exactly two extremal divisorial contractions of the form $f_i\colon X_\C \to Y$, $i=1,2$, where $Y$ is the blow-up of the quadric $Q$ in one conic, see \cite[\S~III.3, p.~80]{Mat95}. If those two contractions are fixed by the Galois action, they descend to $\R$ and induce birational morphisms $h_i\colon X \to W$, where $W$ is a quadric threefold over $\R$ (here the relative Picard rank for $h_1$ is 2). Since $\xr\ne\varnothing$, we deduce that $W(\R)\ne\varnothing$, and therefore $W$ is rational over $\R$. This implies that $X$ is rational. 
\\Assume now that these two extremal contractions $f_1$ and $f_2$ are exchanged by the Galois action. Then we still obtain a birational morphism $g\colon X \to W$ over $\R$, where $W$ is a quadric threefold over $\R$ (in this second case, the relative Picard rank for $g$ is 1). As before, we conclude that $X$ is rational.
\item[{Family \textnumero 3.11}.] $X_\C$ admits exactly three extremal divisorial contractions $f_i\colon X_\C \to Y_i$, where $Y_1$, $Y_2$ and $Y_3$ belong to the families \textnumero 2.25, \textnumero 2.34 and \textnumero 2.35, respectively, see \cite[\S~III.3, p.~81]{Mat95}. This implies that the Galois action fixes all of them, which, as a consequence, descend to $g_i\colon X \to W_i$ over $\R$, where $W_i$ are rational by Proposition~\ref{prop_m=2}. This implies that $X$ is rational.
\item[{Family \textnumero 3.12}.] $X_\C$ admits exactly three extremal divisorial contractions $f_i\colon X_\C \to Y_i$, where $Y_1$, $Y_2$ and $Y_3$ belong to the families \textnumero 2.27, \textnumero 2.33 and \textnumero 2.34, respectively, see \cite[\S~III.3, p.~83]{Mat95}. This implies that the Galois action fixes all of them, which, as a consequence, descend to $g_i\colon X \to W_i$ over $\R$, where $W_i$ are rational by Proposition~\ref{prop_m=2}. This implies that $X$ is rational.
\item[{Family \textnumero 3.13}.] $X_\C$ admits exactly three extremal divisorial contractions $f_i\colon X_\C \to Y$, where $Y$ belongs to the family \textnumero 2.32, see \cite[\S~III.3, p.~84]{Mat95}. Since the Galois group has order two, there exists at least one contraction which is fixed by it, which descends to $g\colon X \to W$ over $\R$, where $W$ is rational by Proposition~\ref{prop_m=2}. Hence $X$ is rational.
\item[{Family \textnumero 3.14}.] $X_\C$ is isomorphic to the blow-up of $\P^3$ in the disjoint union of a plane cubic curve and a point outside the plane and it admits exactly four extremal divisorial contractions, see \cite[\S~III.3, p.~85]{Mat95}. Three of them are of the form $f_i\colon X_\C \to Y_i$, where $Y_1\simeq \P_{\P^2}(\O\oplus \O(1))$, $Y_2\simeq \P_{\P^2}(\O\oplus \O(2))$ and $Y_3$ is the blow-up of $\P^3$ in a plane cubic curve, while the last morphism contracts a divisor to a singular variety. This implies that the Galois action fixes all of them, so they descend to $g_i\colon X \to W_i$ over $\R$, $i=1,2$, where $W_i$ are rational by Proposition~\ref{prop_m=2}. Hence $X$ is rational.
\item[{Family \textnumero 3.15}.] $X_\C$ admits exactly three extremal divisorial contractions $f_i\colon X_\C \to Y_i$, where $Y_1$, $Y_2$ and $Y_3$ belong to the families \textnumero 2.29, \textnumero 2.31 and \textnumero 2.34, respectively, see \cite[\S~III.3, p.~87]{Mat95}. This implies that the Galois action fixes all of them, which, as a consequence, descend to $g_i\colon X \to W_i$ over $\R$, where $W_i$ are rational by Proposition~\ref{prop_m=2}. This implies that $X$ is rational.
\item[{Family \textnumero 3.16}.] $X_\C$ admits exactly three extremal divisorial contractions $f_i\colon X_\C \to Y_i$, where $Y_1$, $Y_2$ and $Y_3$ belong to the families \textnumero 2.27, \textnumero 2.32 and \textnumero 2.35, respectively, see \cite[\S~III.3, p.~88]{Mat95}. This implies that the Galois action fixes all of them, which, as a consequence, descend to $g_i\colon X \to W_i$ over $\R$, where $W_i$ are rational by Proposition~\ref{prop_m=2}. This implies that $X$ is rational.
\item[{Family \textnumero 3.17}.] $X_\C$ admits a $\P^1$-bundle structure $f\colon X_\C \to \P^1\times\P^1$, see \cite[\S~III.3, p.~89]{Mat95}, \cite[Lemma~6.5.1]{BFT23}. This $\P^1$-bundle structure is unique, so it descends to $g\colon X \to W$ over $\R$, where $W$ is a real form of $\P^1\times \P^1$. Since $\xr\ne\varnothing$, we deduce that $W(\R)\ne\varnothing$, and therefore $W$ is rational over $\R$. This implies that $X$ is rational.
\item[{Family \textnumero 3.18}.] $X_\C$ admits exactly three extremal divisorial contractions $f_i\colon X_\C \to Y_i$, where $Y_1$, $Y_2$ and $Y_3$ belong to the families \textnumero 2.29, \textnumero 2.30 and \textnumero 2.33, respectively, see \cite[\S~III.3, p.~90]{Mat95}. This implies that the Galois action fixes all of them, which, as a consequence, descend to $g_i\colon X \to W_i$ over $\R$, where $W_i$ are rational by Proposition~\ref{prop_m=2}. This implies that $X$ is rational.
\item[{Family \textnumero 3.19}.] $X_\C$ is isomorphic to the blow-up of the quadric $Q\subset \P^4_\C$ in two non-colinear points and admits exactly four extremal divisorial contractions. Two of them are of the form $f_i\colon X_\C \to Y$, $i=1,2$, where $Y$ is the blow-up of the quadric $Q$ in one point, see \cite[\S~III.3, p.~91-p.~92]{Mat95}. If those two contractions are fixed by the Galois action, they descend to $\R$ and induce birational morphisms $h_i\colon X \to W$, where $W$ is a quadric threefold over $\R$ (here the relative Picard rank for $h_1$ is 2). Since $\xr\ne\varnothing$, we deduce that $W(\R)\ne\varnothing$, and therefore $W$ is rational over $\R$. This implies that $X$ is rational. 
\\Assume now that these two extremal contractions $f_1$ and $f_2$ are exchanged by the Galois action. Then we still obtain a birational morphism $g\colon X \to W$ over $\R$, where $W$ is a quadric threefold over $\R$ (in this second case, the relative Picard rank for $g$ is 1). As before, we conclude that $X$ is rational. 
\item[{Family \textnumero 3.20}.] $X_\C$ admits exactly three extremal divisorial contractions $f_i\colon X_\C \to Y_i$, where $Y_1$ and $Y_2$ belong to the family \textnumero 2.31 and $Y_3$ belongs to the family \textnumero 2.32, see \cite[\S~III.3, p.~93]{Mat95}. This implies that the Galois action fixes the third one, which, as a consequence, descends to $g_3\colon X \to W_3$ over $\R$, where $W_3$ is rational by Proposition~\ref{prop_m=2}. This implies that $X$ is rational.
\item[{Family \textnumero 3.21}.] $X_\C$ admits a unique extremal divisorial contraction $f\colon X_\C \to \P^1\times \P^2$, see \cite[\S~III.3, p.~94]{Mat95}, so it descends to $g\colon X \to \P^1\times \P^2$ over $\R$. This implies that $X$ is rational.
\item[{Family \textnumero 3.22}.] $X_\C$ admits exactly two extremal divisorial contractions $f_i\colon X_\C \to Y_i$, where $Y_1\simeq \P^1\times \P^2$ and $Y_2\simeq \P_{\P^2}(\O\oplus \O(2))$, see \cite[\S~III.3, p.~96]{Mat95}. Since they cannot be exchanged by the Galois action, they descend to $g_i\colon X \to W_i$, where $W_i$, $i=1,2$, are rational by Proposition~\ref{prop_m=2}. This implies that $X$ is rational.
\item[{Family \textnumero 3.23}.] $X_\C$ admits exactly three extremal divisorial contractions $f_i\colon X_\C \to Y_i$, where $Y_1\simeq \P_{\P^2}(\O\oplus \O(1))$, $Y_2$ belongs to the family \textnumero 2.30 and $Y_2$ belongs to family \textnumero 2.31, see \cite[\S~III.3, p.~97]{Mat95}. Since they cannot be exchanged by the Galois action, they descend to $g_i\colon X \to W_i$  over $\R$, where $W_i$, $i=1,2,3$, are rational by Proposition~\ref{prop_m=2}. This implies that $X$ is rational.
\item[{Family \textnumero 3.24}.] $X_\C$ admits exactly two extremal divisorial contractions $f_i\colon X_\C \to Y_i$, where $Y_1\simeq \P^1\times \P^2$ and $Y_2$ belongs to the family \textnumero 2.32, see \cite[\S~III.3, p.~98]{Mat95}. Since they cannot be exchanged by the Galois action, they descend to $g_i\colon X \to W_i$ over $\R$, where $W_i$, $i=1,2$, are rational by Proposition~\ref{prop_m=2}. This implies that $X$ is rational.
\item[{Family \textnumero 3.25}.] $X_\C$ admits a $\P^1$-bundle structure $f\colon X_\C \to \P^1\times\P^1$, see \cite[\S~III.3, p.~99]{Mat95}, \cite[Lemma~6.5.1]{BFT23}. This $\P^1$-bundle structure is unique, so it descends to $g\colon X \to W$ over $\R$, where $W$ is a real form of $\P^1\times \P^1$. Since $\xr\ne\varnothing$, we deduce that $W(\R)\ne\varnothing$, and therefore $W$ is rational over $\R$. This implies that $X$ is rational.
\item[{Family \textnumero 3.26}.] $X_\C$ admits exactly three extremal divisorial contractions $f_i\colon X_\C \to Y_i$, where $Y_1\simeq \P^1\times \P^2$, $Y_2\simeq \P_{\P^2}(\O\oplus \O(1))$, and $Y_3$ is the blow-up of $\P^3$ at a line, see \cite[\S~III.3, p.~100]{Mat95}. Since they cannot be permuted by the Galois action, they descend to $g_i\colon X \to W_i$ over $\R$, where $W_i$, $i=1,2,3$, are rational by Proposition~\ref{prop_m=2}. This implies that $X$ is rational.
\item[{Family \textnumero 3.27}.] $X_\mathbb{C}=\P^1\times\P^1\times\P^1$. By Lemma~\ref{lem_dp}, $X$ is rational. 
\item[{Family \textnumero 3.28}.] $X_\mathbb{C}=\mathbb{P}^1\times \F_1$.
By Lemma~\ref{lem_dp}, $X$ is rational.
\item[{Family \textnumero 3.29}.] $X_\C$ admits exactly two extremal divisorial contractions $f_i\colon X_\C \to Y_i$, where $Y_1\simeq \P_{\P^2}(\O\oplus \O(1))$ and $Y_2\simeq \P_{\P^2}(\O\oplus \O(2))$, see \cite[\S~III.3, p.~45]{Mat95}. Since they cannot be exchanged by the Galois action, they descend to $g_i\colon X \to W_i$, where $W_i$, $i=1,2$, are rational by Proposition~\ref{prop_m=2}. This implies that $X$ is rational.
\item[{Family \textnumero 3.30}.] $X_\C$ admits a $\P^1$-bundle structure $f\colon X_\C \to \F_1$, see \cite[\S~III.3, p.~103]{Mat95}, \cite[Lemma~6.5.1]{BFT23}. This $\P^1$-bundle structure is unique, so it descends to $g\colon X \to \F_1$ over $\R$ and therefore $X$ is rational.
\item[{Family \textnumero 3.31}.] $X_\C$ admits a $\P^1$-bundle structure $f\colon X_\C \to \P^1\times \P^1$, see \cite[\S~III.3, p.~104]{Mat95}, \cite[Lemma~6.5.1]{BFT23}. This $\P^1$-bundle structure is unique, so it descends to $g\colon X \to W$, where $W$ is a real form of $\P^1\times \P^1$. Since $\xr\ne\varnothing$, we deduce that $W(\R)\ne\varnothing$, and therefore $W$ is rational over $\R$. This implies that $X$ is rational.
\end{itemize}
\endproof

We now study the remaining families with geometric Picard rank three.

\begin{proposition}\label{rang_3_more}
The following holds:
\begin{itemize}
\item[(i)] $s_{3.2}\ge 3$;
\item[(ii)] $s_{3.3}\ge 4$;
\item[(iii)] $s_{3.4}= 3$.
\end{itemize}
\end{proposition}
\proof\begin{itemize}
\item[{Family \textnumero 3.2}.] $X_\C$ is a divisor from $|\mathcal{L}^{\otimes 2}\otimes \mathcal{O}(2,3)|$ on the $\P^2$-bundle $\P(\mathcal{O} \oplus \mathcal{O}(-1,-1)^{\oplus 2})$ over $\P^1\times \P^1$ such that $X_\C\cap Y$ is irreducible, $\mathcal{L}$ is the tautological bundle, and $Y\in |\mathcal{L}|$. 
As explained in \cite[\S~III.3, p.~70-p.~71]{Mat95} and in \cite[Section~5.11]{Chelbook23}, \cite[Lemma 5.11]{ACKM24}, $X_\C$ admits a divisorial contraction onto a non-$\Q$-factorial Fano threefold with one isolated ordinary double point of degree 16: this is precisely the family of threefolds appearing in the proof of Proposition~\ref{prop:X16}. The same proof implies that $s_{3.2}\ge 3$.
\item[{Family \textnumero 3.3}.] $X_\C$ is a divisor on $\P^1\times \P^1 \times \P^2$ of tridegree $(1,1,2)$. Consider the equations
$$x_3\left(-y_0^2-\frac{1}2y_1^2+y_2^2\right)+x_2\left(\frac{1}2y_0^2+y_1^2-y_2^2\right)+\left(x_0+\frac{1}2x_2\right)y_2^2=0$$
$$x_0^2+x_1^2+x_2x_3=0$$
in $\P^3\times \P^2$ with coordinates $([x_0:x_1:x_2:x_3],[y_0:y_1:y_2])$. This defines a (geometrically non-)standard conic bundle over $\P^2$, with discriminant of degree 4, defined by $\left(-y_0^2-\frac{1}2y_1^2+y_2^2\right)\left(\frac{1}2y_0^2+y_1^2-y_2^2\right)-y_2^4=0$. One can check that this curve has four connected components and the corresponding threefold verifies $\#\pi_0(X(\R))=4$. This implies that $s_{3.3}\ge 4$.
\item[{Family \textnumero 3.4}.] $X_\C$ is the blow-up of a theefold $Y$ belonging to family 2.18 in a smooth fiber of the conic bundle structure over $\P^2$. By \cite[\S~III.3, p.~72-p.~73]{Mat95}, the divisorial contraction to $Y$ descends over $\R$, so it is enough to consider a real form in family 2.18 with three connected components, which exists by Proposition~\ref{rang_2_more} and blow up a smooth fiber with nonempty real locus of the conic bundle structure over $\P^2$. This shows that $s_{3.4}= 3$.
\end{itemize}
\endproof

\subsubsection{$\rho(X_\C)=4,5,6$}
\begin{proposition}
\label{prop_m=456}
Let $X$ be a smooth real Fano threefold 
such that $\rho(X_\C)=4,5,6$. Assume $X$ does not belong to family \textnumero 4.1 and $\xr\ne\varnothing$, then $X$ is rational and then $s_X=1$.
\end{proposition}

\proof
First assume that $\rho(X_\C)=4$. There are 13 families of smooth complex Fano threefolds with Picard number~$4$ and any $X$ belonging to any family \textnumero $4.n$ is geometrically rational. We leave aside $n=1$.

If $n\in\{3,4,5,6,8,9,11,12,13\}$, $X_\C$ admits, among its extremal divisorial contractions, a special one $f\colon X_\C \to Y$, where $Y$ is $\P^3$, the quadric $Q\subset \P^4$, the product $\P^1\times\P^1\times\P^1$, 
the blow up of $\P^3$ in two points, the product $\P^1\times\F_1$ 
or the blow up of $Q$ in two non-colinear points. 
From the description of the Mori cone $\overline{\mathrm{NE}}(X_\C)$ given in \cite[\S~III.3, p.106-p.123]{Mat95} and \cite[Section~4]{Mat-erratum}, we deduce that, for every such family, this special contraction descends to a birational morphism $g\colon X \to W$, where $W$ is a real form of $Y$ (the strategy is the same as for the previous results in lower rank). Since $\xr\ne\varnothing$, we have $W(\R)\ne\varnothing$, and therefore $W$ is rational over $\R$ by Proposition~\ref{prop_m=3}, hence $X$ is rational.

If $n=2$, the description of the Mori cone $\overline{\mathrm{NE}}(X_\C)$ given in \cite[\S~III.3, p.107-p.109]{Mat95} shows that $X_\C$ admits two extremal birational contractions $f\colon X_\C \to Y$, where $Y$ belongs to family \textnumero 3.31. Looking at the intersection numebers of the corresponding extremal rays $l_1$ and $l_2$ in the table at \cite[p.108]{Mat95}, we see that the Galois action cannot exchange them. 
Thus $X$ admits an extremal birational contraction $f\colon X \to Y$, defined over $\R$, where $Y$ belongs to family \textnumero 3.31. Since $\xr\ne\varnothing$, we have $Y(\R)\ne\varnothing$, and therefore $Y$ is rational over $\R$ by Proposition~\ref{prop_m=3}.

If $n=7$, the description of the Mori cone $\overline{\mathrm{NE}}(X_\C)$ given in \cite[\S~III.3, p.116-p.117]{Mat95} shows that $X$ admits a real birational map $f\colon X \to Y$, where $Y$ belongs to family \textnumero 2.32. Since $\xr\ne\varnothing$, we have $Y(\R)\ne\varnothing$, and therefore $Y$ is rational over $\R$ by Proposition~\ref{prop_m=2}. 

If $n=10$, then $X_\mathbb{C}=\mathbb{P}^1\times S_7$.
By Lemma~\ref{lem_dp}, $X$ is rational.

\medskip
Now assume $\rho(X_\C)=5$ or $6$.
There are 4 families of smooth complex Fano threefolds with Picard number $5$ or $6$ and any $X$ belonging to one of these families is geometrically rational.
Assume that $X_\C$ belongs to family \textnumero 5.1, then it admits three extremal birational contractions $f\colon X_\C \to Y$, where $Y$ belongs to family \textnumero 4.12.
We deduce from the description of the Mori cone $\overline{\mathrm{NE}}(X_\C)$ given in \cite[\S~III.3, p.124-p.125]{Mat95} that at least one of them is defined over $\R$. Since $\xr\ne\varnothing$, then $Y(\R)\ne\varnothing$, thus $Y$ is rational over $\R$ by the previous part of the proof.

Assume that $X_\C$ belongs to family \textnumero 5.2, then, among its extremal birational contractions, $X_\C$ admits a unique one $f\colon X_\C \to Y$, whose target space $Y$ belongs to family \textnumero 4.12 (see the description of the Mori cone $\overline{\mathrm{NE}}(X_\C)$ given in \cite[\S~III.3, p.125-p.126]{Mat95}). We deduce 
that $f$ descends over $\R$. Since $\xr\ne\varnothing$, then $Y(\R)\ne\varnothing$, thus $Y$ is rational over $\R$ by the previous part of the proof. 

Assume that $X_\C$ belongs to family \textnumero 5.3 or \textnumero 6.1, then $X_\mathbb{C}=\mathbb{P}^1\times S_{11-m}$ and by Lemma~\ref{lem_dp}, $X$ is rational.
\endproof

\begin{proposition}
\label{rang_456_more}
The following holds: $s_{4.1}=2$.
\end{proposition}

\proof
By \cite[Example 4.3]{CTZ24-pre}, there exists a smooth real Fano threefold $X$ belonging to family \textnumero 4.1 whose real locus is disconnected, hence $s_{4.1}\geq2$. Conversely let $X$ be a smooth real Fano threefold belonging to family \textnumero 4.1. Then $X$ is a real form of a degree $(1,1,1,1)$ divisor in $\P^1\times\P^1\times\P^1\times\P^1$. In fact, any such form is a degree $(1,1,1,1)$ divisor in a real form $Y$ of $\P^1\times\P^1\times\P^1\times\P^1$. Exactly three of these forms have real points and $X$ can be non rational  iff $Y=Q\times Q$ where $Q$ is the smooth quadric surface whose real locus is a sphere. 
In this case, we have $\rho(X)=2$ and $\rho(X_\C)=4$ and by Proposition~\ref{prop.st.gm}, we get $s_{4.1}\leq 2$ taking into account that $h^{1,2}(X)=1$. 
\endproof

\subsubsection{$\rho(X_\C)\ge 7$}

\begin{proposition}
\label{proposition:78910}
Let $X$ be a smooth real Fano threefold such that $\rho(X_\C)\ge 7$, then $s_X>1$. More precisely, $s_{7.1}=s_{8.1}=2$, $s_{9.1}=4$, $s_{10.1}=5$.
\end{proposition}
\proof
There are 4 families of smooth complex Fano threefolds with Picard number $m\ge 7$ and if $X$ belongs to one of these families, $X_\mathbb{C}=\mathbb{P}^1\times S_{11-m}$ thus Lemma~\ref{lem_dp} applies and we get the result by Example~\ref{example:sx}\ref{item:dp}.
\endproof

\subsubsection{Final remarks}

\proof[Proof of Theorem~\ref{thm_newmain}]
It follows directly from Propositions~\ref{prop_m=1}, \ref{rang_1_more}, \ref{prop_m=2}, \ref{rang_2_more}, \ref{prop_m=3}, \ref{rang_3_more}, \ref{prop_m=456}, \ref{rang_456_more} and \ref{proposition:78910}.
\endproof

\begin{remark}
If $m >6$, we have shown above that $s_{m.n}>1$ but $X$ is rational if and only if its real locus is nonempty and connected by Lemma~\ref{lem_dp} and rational examples exist. Thus, the converse of Theorem~\ref{thm_newmain} is not true.
\end{remark}

\begin{remark}
We do not know whether, in Theorem~\ref{thm_newmain}, we can relax the Fano hypothesis. Indeed in \cite{BP24-pre} one can find a family of real conic bundles  
over rational surfaces
whose elements are irrational with connected real locus. 
Namely, let $U:=\{x^2+y^2=f(u,v)\}\subset \A^4$ where $f\in\R[u,v]$ is a polynomial of even degree $d\geq 12$ such that the closure $C$ of $\{f = 0\}$ in $\P^2$ is a nodal rational curve. Blowing up the nodes of the curve $C$, we get a rational surface $S$ and a real conic bundle $X\to S$ which is a smooth projective models of $U$. Then if the real locus $X(\R)$ is nonempty, it is connected but at the same time it is irrational by \cite[Theorem~4.5]{BP24-pre}.

To prove that the statement of Theorem~\ref{thm_newmain} holds true for this family of examples, one needs to provide a member $X$ of this family of conic bundles with non-trivial automorphism group and a twisted real form $Y$ of $X$ with non-connected real locus $Y(\R)$.
\end{remark}

\begin{question}
Let $n\ge 3$ be an integer and $X$ be a smooth geometrically rational real $n$-dimensional variety with $\xr\ne\varnothing$ and assume that $s_X=1$. Is $X$ rational?
\end{question}

The answer to the same question is positive for $n=1$ (trivial), and for $n=2$ by Comessatti's Theorem.


\section{Recap Table}\label{sec_table}

The following table summarizes the results for 
smooth geometrically rational real Fano threefold for which some deformation of $X_\C$ admits a real form whose real locus has at least two connected components. 
In particular, columns $s_{m.n}\geq$ and $s_{m.n}\leq$ recap the bounds collected in this paper. We deduce the lower bound producing examples, while the upper bound is the minimum between the Smith-Thom and Borel-Swan inequalities or via a classification. 

We use the following notation: the numbering of families \textnumero $m.n$ is the one in \cite{fanography}, and for $X$ belonging to the family \textnumero $m.n$, $\iota$ denotes the Fano index of $X_\C$, $d=(-K_X^3)/\iota^3$ its degree, $h^{1,2}:=h^{1,2}(X_\C)$ the corresponding Hodge number, the definition of the maximal number of connected components $s_{m.n}=s_X$ can be found in Definition~\ref{dfn.smn}, ``$\exists$ IC'' stands for the existence of an irrational real Fano threefold in the family with nonempty connected real locus. The right column gives a description of $X_\C$. Some of those descriptions remain true over $\R$, assuming that $X(\R)$ is non-empty, namely \textnumero $1.14$, $2.10$, $2.18$, $3.4$, $7.1$, $8.1$, $9.1$ and $10.1$.

{
\begin{center}
\phantomsection\refstepcounter{tables}
\renewcommand{\arraystretch}{1.2}
\begin{longtable}{@{}l|c|c|c|c|c|c|l@{}}
\toprule
\textnumero $m.n$ & $\iota$	& $d$ & $h^{1,2}$	& $s_{m.n}\geq$ & $s_{m.n}\leq$ & $\exists$ IC
& Description of $X_\C$
 \\ 
\midrule
\endfirsthead
\toprule
\textnumero $m.n$ & $\iota$	& $d$ & $h^{1,2}$	& $s_{m.n}\geq$ & $s_{m.n}\leq$ & $\exists$ IC
& Description of $X_\C$
 \\ 
\midrule
\endhead
\midrule \multicolumn{8}{r}{\textit{Continued on next page}} \\
\endfoot
\endlastfoot
1.8 	& 1	& 16		& 3	& 3    & 5    & ?     & $X_{16}\subset \P^{10}$, section of $\LGr(2,5)\subset \P^{13}$ by 
\\& & & & & & & a linear subspace of codimension 3
\\ 
1.14	& 2	& 4		& 2	& 2     & 2   & yes  & $V_4\subset \P^5$, smooth complete intersection 
\\& & & & & & & of two quadrics
\\
2.10	& 1	& 16		& 3	& 2     & 2	   & yes  & $\Bl_CV_4$, where $C$ is an elliptic curve
\\
2.12	& 1	& 20		& 3	& 2     & 4    & yes & $X_{(3,3)}\subset\mathbb{P}^3\times \mathbb{P}^3$ smooth intersection of \\ & & & & & & & three $(1,1)$-divisors
\\
2.18	& 1	& 24		& 2	& 3     & 3     & yes &
$X_{24}$, double cover of $\mathbb{P}^1\times \mathbb{P}^2$ branched 
\\& & & & & & & at a $(2,2)$-divisor
\\
3.2	& 1	& 14  	& 3	& 3     & 7     & ?    & Divisor in $|\mathcal{L}^{\otimes 2}\otimes \mathcal{O}(2,3)|$ on the  
\\ & & & & & & & $\P^2$-bundle $\P(\mathcal{O} \oplus \mathcal{O}(-1,-1)^{\oplus 2})$ over 
\\ & & & & & & & $\P^1\times \P^1$, with $X_\C\cap Y$ is irreducible, 
\\ & & & & & & & $\mathcal{L}$ tautological bundle, $Y\in |\mathcal{L}|$
\\
3.3	& 1	& 18  	& 3	& 4     & 5     & ?    & $(1,1,2)$-divisor in $\mathbb{P}^1\times\mathbb{P}^1\times\mathbb{P}^2$
\\
3.4 	& 1	& 18  	& 2	& 3     & 3     & yes & $\Bl_FX_{24}$, where $F$ is a smooth fiber of 
\\ & & & & & & & the composition of the projection to 
\\ & & & & & & & $\mathbb{P}^1\times\mathbb{P}^2$ with the projection to $\mathbb{P}^2$
\\
4.1	& 	& 24  	& 1    & 2     & 2     & ?    & $(1,1,1,1)$-divisor in $(\mathbb{P}^1)^4$
\\
7.1	& 	& 24  	& 0    & 2     & 2      & no & $\mathbb{P}^1\times S_4$
\\
8.1	& 	& 18  	& 0    & 2     & 2      & no & $\mathbb{P}^1\times S_3$
\\
9.1	& 	& 12  	& 0    & 4     & 4      & no & $\mathbb{P}^1\times S_2$
\\
10.1	& 	& 6  		& 0    & 5     & 5      & no & $\mathbb{P}^1\times S_1$
\\[1mm]
\bottomrule
\multicolumn{8}{c}{\textrm{Table~\ltables{tab.recap}. Real geometrically rational Fano threefolds with $s_{m.n}>1$}}\\[1mm]
\end{longtable}
\end{center}
}

\bibliographystyle{alpha}
\bibliography{biblio}

\newcommand{\etalchar}[1]{$^{#1}$}
\begin{thebibliography}{CTSSD87}

\bibitem[AB92]{AB92}
Alberto Alzati and Marina Bertolini.
\newblock On the rationality of {F}ano {$3$}-folds with {$B_2\geq 2$}.
\newblock {\em Matematiche (Catania)}, 47(1):63--74 (1993), 1992.

\bibitem[ACC{\etalchar{+}}23]{Chelbook23}
Carolina Araujo, Ana-Maria Castravet, Ivan Cheltsov, Kento Fujita, Anne-Sophie
  Kaloghiros, Jesus Martinez-Garcia, Constantin Shramov, Hendrik S\"u\ss, and
  Nivedita Viswanathan.
\newblock {\em The {C}alabi problem for {F}ano threefolds}, volume 485 of {\em
  London Mathematical Society Lecture Note Series}.
\newblock Cambridge University Press, Cambridge, 2023.

\bibitem[ACKM24]{ACKM24}
Hamid Abban, Ivan Cheltsov, Takashi Kishimoto, and Fr{\'e}d{\'e}ric Mangolte.
\newblock {K}-stability of pointless del {P}ezzo surfaces and {F}ano 3-folds.
\newblock {\em preprint, arXiv:2411.00767}, 2024.

\bibitem[ACKM25]{ackm3}
Hamid Abban, Ivan Cheltsov, Takashi Kishimoto, and Fr{\'e}d{\'e}ric Mangolte.
\newblock {S}mooth {F}ano 3-folds satisfying {C}ondition $(\mathbf{A})$.
\newblock Preprint, \texttt{arXiv:2505.13684}, 2025.

\bibitem[Bel15]{fanography}
P.~Belmans.
\newblock Fanography.
\newblock \url{https://fanography.info}, 2015.

\bibitem[BFT23]{BFT23}
J\'er\'emy Blanc, Andrea Fanelli, and Ronan Terpereau.
\newblock Connected algebraic groups acting on three-dimensional {M}ori
  fibrations.
\newblock {\em Int. Math. Res. Not. IMRN}, pages 1572--1689, 2023.

\bibitem[BP24]{BP24-pre}
Olivier Benoist and Alena Pirutka.
\newblock On the rationality of some real threefolds.
\newblock Preprint, \texttt{2412.13624}, 2024.

\bibitem[BW20]{BW20}
Olivier Benoist and Olivier Wittenberg.
\newblock The {C}lemens-{G}riffiths method over non-closed fields.
\newblock {\em Algebr. Geom.}, 7(6):696--721, 2020.

\bibitem[BW23]{BW23}
Olivier Benoist and Olivier Wittenberg.
\newblock Intermediate {J}acobians and rationality over arbitrary fields.
\newblock {\em Ann. Sci. \'Ec. Norm. Sup\'er. (4)}, 56(4):1029--1084, 2023.

\bibitem[CFST16]{CFST16}
Giulio Codogni, Andrea Fanelli, Roberto Svaldi, and Luca Tasin.
\newblock Fano varieties in {M}ori fibre spaces.
\newblock {\em Int. Math. Res. Not. IMRN}, pages 2026--2067, 2016.

\bibitem[CFST18]{CFST18}
Giulio Codogni, Andrea Fanelli, Roberto Svaldi, and Luca Tasin.
\newblock A note on the fibres of {M}ori fibre spaces.
\newblock {\em Eur. J. Math.}, 4(3):859--878, 2018.

\bibitem[CM08]{cm1}
Fabrizio Catanese and Fr{\'e}d{\'e}ric Mangolte.
\newblock Real singular del {P}ezzo surfaces and 3-folds fibred by rational
  curves. {I}.
\newblock {\em Michigan Math. J.}, 56(2):357--373, 2008.

\bibitem[CM09]{cm2}
Fabrizio Catanese and Fr{\'e}d{\'e}ric Mangolte.
\newblock Real singular del {P}ezzo surfaces and 3-folds fibred by rational
  curves. {II}.
\newblock {\em Ann. Sci. \'Ec. Norm. Sup\'er. (4)}, 42(4):531--557, 2009.

\bibitem[Com12]{Co12}
Annibale Comessatti.
\newblock Fondamenti per la geometria sopra le suerficie razionali dal punto di
  vista reale.
\newblock {\em Math. Ann.}, 73(1):1--72, 1912.

\bibitem[CTP24]{CP24}
Jean-Louis Colliot-Thélène and Alena Pirutka.
\newblock Certaines fibrations en surfaces quadriques réelles.
\newblock {\em preprint, arXiv:2406.00463}, 2024.

\bibitem[CTSSD87]{CTSSD87}
Jean-Louis Colliot-Th\'{e}l\`ene, Jean-Jacques Sansuc, and Peter
  Swinnerton-Dyer.
\newblock Intersections of two quadrics and {C}h\^{a}telet surfaces. {II}.
\newblock {\em J. Reine Angew. Math.}, 374:72--168, 1987.

\bibitem[CTZ24]{CTZ24-pre}
Ivan Cheltsov, Yuri Tschinkel, and Zhang Zhijia.
\newblock Rationality of singular cubic threefolds over $\mathbb{R}$.
\newblock {\em preprint, arXiv:2411.14379}, 2024.

\bibitem[Deb20]{De20}
Olivier Debarre.
\newblock {G}ushel-{M}ukai varieties.
\newblock preprint, arXiv:2001.03485, 2020.

\bibitem[DK19]{DK19}
Adrien Dubouloz and Takashi Kishimoto.
\newblock Cylindres dans les fibrations de {M}ori: formes du volume quintique
  de del {P}ezzo.
\newblock {\em Ann. Inst. Fourier (Grenoble)}, 69(6):2377--2393, 2019.

\bibitem[EGH00]{EGH}
Y.~Eliashberg, A.~Givental, and H.~Hofer.
\newblock Introduction to symplectic field theory.
\newblock {\em Geom. Funct. Anal.}, Special Volume, Part II:560--673, 2000.
\newblock GAFA 2000 (Tel Aviv, 1999).

\bibitem[FJ24]{FJ24}
Sarah Frei and Lena Ji.
\newblock A threefold violating a local-to-global principle for rationality.
\newblock {\em Res. Number Theory}, 10(2):Paper No. 39, 9, 2024.

\bibitem[FJS{\etalchar{+}}24a]{FJSVV2}
Sarah Frei, Lena Ji, Soumya Sankar, Bianca Viray, and Isabel Vogt.
\newblock Conic bundle threefolds differing by a constant brauer class and
  connections to rationality.
\newblock {\em preprint, arXiv:2406.13510}, 2024.

\bibitem[FJS{\etalchar{+}}24b]{FJSVV24}
Sarah Frei, Lena Ji, Soumya Sankar, Bianca Viray, and Isabel Vogt.
\newblock Curve classes on conic bundle threefolds and applications to
  rationality.
\newblock {\em Algebr. Geom.}, 11(3):421--459, 2024.

\bibitem[GS17]{GS17}
Philippe Gille and Tam\'as Szamuely.
\newblock {\em Central simple algebras and {G}alois cohomology}, volume 165 of
  {\em Cambridge Studies in Advanced Mathematics}.
\newblock Cambridge University Press, Cambridge, second edition, 2017.

\bibitem[Har76]{Kh76}
V.~M. Harlamov.
\newblock Topological types of nonsingular surfaces of degree {$4$} in {${\bf
  R}P\sp{3}$}.
\newblock {\em Funkcional. Anal. i Prilo\v zen.}, 10(4):55--68, 1976.

\bibitem[HM05a]{hm2}
Johannes Huisman and Fr{\'e}d{\'e}ric Mangolte.
\newblock Every connected sum of lens spaces is a real component of a uniruled
  algebraic variety.
\newblock {\em Ann. Inst. Fourier (Grenoble)}, 55(7):2475--2487, 2005.

\bibitem[HM05b]{hm1}
Johannes Huisman and Fr{\'e}d{\'e}ric Mangolte.
\newblock Every orientable {S}eifert 3-manifold is a real component of a
  uniruled algebraic variety.
\newblock {\em Topology}, 44(1):63--71, 2005.

\bibitem[HT21a]{HT21}
Brendan Hassett and Yuri Tschinkel.
\newblock Cycle class maps and birational invariants.
\newblock {\em Comm. Pure Appl. Math.}, 74(12):2675--2698, 2021.

\bibitem[HT21b]{HT21b}
Brendan Hassett and Yuri Tschinkel.
\newblock Rationality of complete intersections of two quadrics over nonclosed
  fields.
\newblock {\em Enseign. Math.}, 67(1-2):1--44, 2021.
\newblock With an appendix by Jean-Louis Colliot-Th\'{e}l\`ene.

\bibitem[HT22]{HT19}
Brendan Hassett and Yuri Tschinkel.
\newblock Rationality of {F}ano threefolds of degree 18 over nonclosed fields.
\newblock {\em preprint, arXiv:1910.13816v1}, 2022.

\bibitem[IP99]{IP99}
V.~A. Iskovskikh and Yu.~G. Prokhorov.
\newblock Fano varieties.
\newblock In {\em Algebraic geometry, {V}}, volume~47 of {\em Encyclopaedia
  Math. Sci.}, pages 1--247. Springer, Berlin, 1999.

\bibitem[Isk79]{I79}
V.~A. Iskovskih.
\newblock Anticanonical models of three-dimensional algebraic varieties.
\newblock In {\em Current problems in mathematics, {V}ol. 12 ({R}ussian)},
  pages 59--157, 239 (loose errata). VINITI, Moscow, 1979.

\bibitem[JJ24]{JJ23}
Lena Ji and Mattie Ji.
\newblock Rationality of real conic bundles with quartic discriminant curve.
\newblock {\em Int. Math. Res. Not. IMRN}, 2024(1):115--151, 2024.

\bibitem[JR11]{JR11}
Priska Jahnke and Ivo Radloff.
\newblock Terminal {F}ano threefolds and their smoothings.
\newblock {\em Math. Z.}, 269(3-4):1129--1136, 2011.

\bibitem[KI96]{IK96}
V.~Kharlamov and I.~Itenberg.
\newblock Towards the maximal number of components of a nonsingular surface of
  degree {$5$} in {${\bf R}{\rm P}^3$}.
\newblock In {\em Topology of real algebraic varieties and related topics},
  volume 173 of {\em Amer. Math. Soc. Transl. Ser. 2}, pages 111--118. Amer.
  Math. Soc., Providence, RI, 1996.

\bibitem[Kol98a]{K98}
J\'{a}nos Koll\'{a}r.
\newblock The {N}ash conjecture for threefolds.
\newblock {\em Electron. Res. Announc. Amer. Math. Soc.}, 4:63--73, 1998.

\bibitem[Kol98b]{KoI}
J{\'a}nos Koll{\'a}r.
\newblock Real algebraic threefolds. {I}. {T}erminal singularities.
\newblock {\em Collect. Math.}, 49(2-3):335--360, 1998.
\newblock Dedicated to the memory of Fernando Serrano.

\bibitem[Kol99a]{KoII}
J{\'a}nos Koll{\'a}r.
\newblock Real algebraic threefolds. {II}. {M}inimal model program.
\newblock {\em J. Amer. Math. Soc.}, 12(1):33--83, 1999.

\bibitem[Kol99b]{KoIII}
J{\'a}nos Koll{\'a}r.
\newblock Real algebraic threefolds. {III}. {C}onic bundles.
\newblock {\em J. Math. Sci. (New York)}, 94(1):996--1020, 1999.
\newblock Algebraic geometry, 9.

\bibitem[Kol00]{KoIV}
J{\'a}nos Koll{\'a}r.
\newblock Real algebraic threefolds. {IV}. {D}el {P}ezzo fibrations.
\newblock In {\em Complex analysis and algebraic geometry}, pages 317--346. de
  Gruyter, Berlin, 2000.

\bibitem[Kol01]{K01}
J\'{a}nos Koll\'{a}r.
\newblock The topology of real algebraic varieties.
\newblock In {\em Current developments in mathematics, 2000}, pages 197--231.
  Int. Press, Somerville, MA, 2001.

\bibitem[KP23]{KP23}
Alexander Kuznetsov and Yuri Prokhorov.
\newblock Rationality of {F}ano threefolds over non-closed fields.
\newblock {\em Amer. J. Math.}, 145(2):335--411, 2023.

\bibitem[KP24]{KP24}
Alexander Kuznetsov and Yuri Prokhorov.
\newblock Rationality over nonclosed fields of {F}ano threefolds with higher
  geometric {P}icard rank.
\newblock {\em J. Inst. Math. Jussieu}, 23(1):207--247, 2024.

\bibitem[Kra09]{Kra09}
V.~A. Krasnov.
\newblock On the topological classification of real three-dimensional cubics.
\newblock {\em Mat. Zametki}, 85(6):886--893, 2009.

\bibitem[Kra18]{K18}
V.~A. Krasnov.
\newblock On the intersection of two real quadrics.
\newblock {\em Izv. Ross. Akad. Nauk Ser. Mat.}, 82(1):97--150, 2018.

\bibitem[KS04]{KS04}
J{\'a}nos Koll{\'a}r and Frank-Olaf Schreyer.
\newblock Real {F}ano 3-folds of type {$V_{22}$}.
\newblock In {\em The {F}ano {C}onference}, pages 515--531. Univ. Torino,
  Turin, 2004.

\bibitem[Man14]{M14}
Fr{\'e}d{\'e}ric Mangolte.
\newblock Topologie des variétés algébriques réelles de dimension~3.
\newblock {\em Gaz. Math.}, 139:5--34, 2014.

\bibitem[Man17]{M17}
Fr\'ed\'eric Mangolte.
\newblock {\em Vari\'et\'es alg\'ebriques r\'eelles}, volume~24 of {\em Cours
  Sp\'ecialis\'es [Specialized Courses]}.
\newblock Soci\'et\'e{} Math\'ematique de France, Paris, 2017.

\bibitem[Man20]{M20}
Fr\'ed\'eric Mangolte.
\newblock {\em Real algebraic varieties}.
\newblock Springer Monographs in Mathematics. Springer, Cham, 2020.
\newblock Translated from the 2017 French original by Catriona Maclean.

\bibitem[Mat95]{Mat95}
Kenji Matsuki.
\newblock Weyl groups and birational transformations among minimal models.
\newblock {\em Mem. Amer. Math. Soc.}, 116(557):vi+133, 1995.

\bibitem[Mat23]{Mat-erratum}
Kenji Matsuki.
\newblock {A}ddendum/{E}rratum to the paper "{W}eyl groups and {B}irational
  transformations among minimal models".
\newblock preprint, arXiv:2401.13431, 2023.

\bibitem[MM86]{MM86}
Shigefumi Mori and Shigeru Mukai.
\newblock Classification of {F}ano {$3$}-folds with {$B_2\geq 2$}. {I}.
\newblock In {\em Algebraic and topological theories ({K}inosaki, 1984)}, pages
  496--545. Kinokuniya, Tokyo, 1986.

\bibitem[MM03]{MM03}
Shigefumi Mori and Shigeru Mukai.
\newblock Erratum: ``{C}lassification of {F}ano 3-folds with {$B_2\geq 2$}''
  [{M}anuscripta {M}ath. {\bf 36} (1981/82), no. 2, 147--162].
\newblock {\em Manuscripta Math.}, 110(3):407, 2003.

\bibitem[MW12]{mw1}
Fr{\'e}d{\'e}ric Mangolte and Jean-Yves Welschinger.
\newblock Do uniruled six-manifolds contain {S}ol {L}agrangian submanifolds?
\newblock {\em Int. Math. Res. Not. IMRN}, 2012(7):1569--1602, 2012.

\bibitem[Nam97]{Na97}
Yoshinori Namikawa.
\newblock Smoothing {F}ano {$3$}-folds.
\newblock {\em J. Algebraic Geom.}, 6(2):307--324, 1997.

\bibitem[Ore01]{Or01}
S.~Yu. Orevkov.
\newblock Real quintic surface with 23 components.
\newblock {\em C. R. Acad. Sci. Paris S\'er. I Math.}, 333(2):115--118, 2001.

\bibitem[Pro13]{P13}
Yuri Prokhorov.
\newblock {$G$}-{F}ano threefolds, {II}.
\newblock {\em Adv. Geom.}, 13(3):419--434, 2013.

\bibitem[Sil89]{S89}
Robert Silhol.
\newblock {\em Real algebraic surfaces}, volume 1392 of {\em Lecture Notes in
  Mathematics}.
\newblock Springer-Verlag, Berlin, 1989.

\bibitem[Tak89]{T89}
Kiyohiko Takeuchi.
\newblock Some birational maps of {F}ano {$3$}-folds.
\newblock {\em Compositio Math.}, 71(3):265--283, 1989.

\bibitem[Tak22]{T22}
Kiyohiko Takeuchi.
\newblock Weak {F}ano threefolds with del {P}ezzo fibration.
\newblock {\em Eur. J. Math.}, 8(3):1225--1290, 2022.

\bibitem[Vit99]{V99}
Claude Viterbo.
\newblock Symplectic real algebraic geometry.
\newblock Unpublished, 1999.

\end{thebibliography}

\end{document}